\input btxmac
\baselineskip = 1.1\baselineskip 
\parskip = 2pt


\font\titlefont=cmssdc10 scaled \magstep1
\font\authorfont=cmcsc10 
\font\headfont=cmcsc10
\font\thanksfont=cmr8
\font\small=cmr9
\font\smallbf=cmbx9
\font\smallss=cmss8
\font\smallssi=cmssi9

\font\smalltt=cmtt9
\newcount\sectno
\newcount\parno
\newcount\myeqno

\def\sectnum{\the\sectno}
\def\parnum{\sectnum .{\the\parno}}

\outer\def\mysection#1{\vskip1.5pc plus.75pc \penalty-250
                 \parno=0  \advance\sectno by1
                 \noindent {\headfont \sectnum.\hskip 1pc #1}
                 \nobreak \bigskip \noindent}

\outer\def\subsection#1{\medbreak \noindent {\bf #1.} \nobreak \medskip \noindent}

\outer\def\newpar#1. #2\par {\advance \parno by1 \medbreak \noindent
     {\bf  #1 \parnum. \enspace}#2 \par
     \ifdim\lastskip<\medskipamount \removelastskip\penalty55\medskip\fi}

\outer\def\newclaim#1. #2\par {\advance \parno by1 \medbreak \noindent
     {\bf  #1 \parnum. \enspace}{\sl #2}\par
     \ifdim\lastskip<\medskipamount \removelastskip\penalty55\medskip\fi}

\outer\def\justpar#1. #2\par {\medbreak \noindent
     {\bf  #1. \enspace}#2 \par
     \ifdim\lastskip<\medskipamount \removelastskip\penalty55\medskip\fi}

\outer\def\justclaim#1. #2\par {\medbreak \noindent
     {\bf  #1. \enspace}{\sl #2}\par
     \ifdim\lastskip<\medskipamount \removelastskip\penalty55\medskip\fi}

\def\supaccent#1#2{\leavevmode
 \setbox0=\hbox{#2}\dimen0=\ht0
 \advance\dimen0 by-1ex \rlap{\raise.69\dimen0
 \hbox to 1.2\wd0{\hss\char#1\hss}}\box0{}}

\def\frames{{\supaccent{'27}{$\cal M$}}_{2n \times n}}

\def\Proof{\noindent {\sl Proof.}\enspace}

\def\qedmark{\hbox{\vrule height 4pt width 3pt}}
\def\qedskip{\vrule height 4pt width 0pt depth 1pc}
\def\qed{\penalty 1000\quad\penalty 1000{\qedmark\qedskip}}

\def\R{{\rm I\! R}}
\def\RP{{\rm I\!R I\!P}}
\def\frac#1#2{{#1 \over #2}}

\def\join#1#2{(#1 | #2)}
\def\hom#1#2{\hbox{\rm Hom}(#1,#2)}

\def\tr{\hbox{\rm tr}}
\def\trans{{\scriptscriptstyle T}}
\def\grassmann{\hbox{\rm G}_n}
\def\lagrangian{\Lambda_n}
\def\lgroup#1{\hbox{\rm GL}(#1)}
\def\lalgebra{{\cal G}l(2n)}



\def\kfan#1#2{\hbox{\rm J}_{f}^{#1}(\R ; #2)}

\def\kmaps#1{\hbox{\rm J}^{#1}(\R ;\lgroup{n})}





\def\bx{\hbox{\bf x}}

\def\bv{\hbox{\bf v}}
\def\bw{\hbox{\bf w}}
\def\bp{\hbox{\bf p}}
\def\bT{\hbox{\bf T}}
\def\bS{\hbox{\bf S}}
\def\bF{\hbox{\bf F}}
\def\bG{\hbox{\bf G}}
\def\bA{\hbox{\bf A}}
\def\bB{\hbox{\bf B}}
\def\bX{\hbox{\bf X}}
\def\bI{\hbox{\bf I}}

\def\bJ{\hbox{\bf J}}
\def\bN{\hbox{\bf N}}
\def\bK{\hbox{\bf K}}
\def\bP{\hbox{\bf P}}
\def\cA{{\cal A}}
\def\cB{{\cal B}}

\def\cD{{\cal D}}
\def\cH{{\cal H}}

\def\bzero{{\bf 0}}

\def\mzero{\hbox{\bf {\rm O}}}
\def\point{\item{$\bullet$}}
\def\exp{\hbox{\rm exp}}

\centerline{\titlefont Geometric Invariants of Fanning Curves}
\bigskip
\centerline{{\authorfont J.~C. \'Alvarez Paiva}\footnote*{\thanksfont Partially
funded by FAPESP grant $\hbox{\rm N}^{\circ}$ 2004/01509-0}
{\authorfont and  C.~E. Dur\'an}}
\vskip .5cm

\bigskip
\medskip
\hfil {\hsize = 10cm \vbox{\noindent {\smallbf Abstract.} \small
We study the geometry of an important class of generic curves in the Grassmann
manifolds of $n$-dimensional subspaces and Lagrangian subspaces of $\R^{2n}$
under the action of the linear and linear symplectic groups.
}} \hfil

\bigskip
\hfil {\hsize = 7cm \vbox
{\frenchspacing \smallssi On pr\'ef\`ere calculer plut\^ot que voir :
c'est moins p\'enible et plus convaincant.}\nonfrenchspacing}
\rightline{\smallss --- REN\'E THOM}

\bigskip

\mysection{Introduction}
Curves on Grassmann manifolds often appear in geometry and dynamics through the
following construction: let $\pi : E \rightarrow M$ be a fiber bundle over a
manifold $M$ and let $\phi: \R \times E \rightarrow E$ be a flow.
If $VE \subset TE$ denotes the vertical subbundle (i.e., the kernel of $D\pi$) and
$e$ is a point in $E$, $t \mapsto D\phi_{-t}(V_{\phi_t(e)}E)$ is a curve of
subspaces of $T_e E$. For example, in Riemannian and Finsler geometry
$E = {\cal T}M$ is the punctured tangent bundle (i.e., without the zero section)
of a manifold $M$ and $\phi$ is the geodesic flow.

In this paper we introduce a novel, comprehensive approach to the geometry of
the class of curves in the Lagrangian Grassmannian $\lagrangian$ and the
Grassmannian $\grassmann$ of $n$-dimensional subspaces of $\R^{2n}$ that
arise in the study of semi-sprays and Lagrangian flows. Besides its simplicity
and elegance, this new approach uncovers the geometry behind the formalisms of
J.~Klein and A.~Voutier~\cite{Klein-Voutier:1968}, J.~Grifone~\cite{Grifone:1972a},
and P.~Foulon~\cite{Foulon:1986} for the construction of connections and curvature
for Finsler metrics, Lagrangian systems, and semi-sprays.

The relationship between geodesic flows, Sturm systems, and the projective
differential geometry of curves in the Grassmannians $\grassmann$ and $\lagrangian$
is classical and has been developed---for the most part implicitly---from
different viewpoints. The explicit development of the subject seems to have been
started by S.~Ahdout~(\cite{Ahdout:1989}) who studied it in the context of geodesic
flows of Riemannian metrics and billiard maps. Later it was independently studied by
V.~Ovsienko (\cite{Ovsienko:1989}, \cite{Ovsienko:1990}, and \cite{Ovsienko:1993})
in the context of Sturm systems, and by A.~Agrachev together with R.~Gamkrelidze,
N.~Chtcherbakova, and I.~Zelenko in the context of control theory (see
\cite{Agrachev:1998}, \cite{Agrachev-etal:2004}, \cite{Agrachev-Gamkrelidze:1997},
\cite{Agrachev-Gamkrelidze:2004}, \cite{Agrachev-Zelenko:2001}, \cite{Agrachev-Zelenko:2002a},
and \cite{Agrachev-Zelenko:2002b}).
Our approach encompasses, unifies, and simplifies those of our predecessors,
while at the same time it departs from them in that we not only study
the projective geometry of curves on $\grassmann$ and $\lagrangian$, we study
the genesis and geometry of the invariants themselves.

Using the canonical identification between the tangent space of the
Grassmannian of $n$-dimensional subspaces of $\R^{2n}$ at an $n$-dimensional
subspace $\ell$ and the space  of linear maps from $\ell$ to the quotient
space $\R^{2n}/\ell$, the class of curves we are interested in can be defined as
follows:

\newpar Definition. A smooth curve $\ell(t)$ of $n$-dimensional subspaces of
$\R^{2n}$ is said to be {\sl fanning\/} if at each time $t$ the tangent
vector $\dot{\ell}(t)$ is an invertible linear map from $\ell(t)$ to the quotient
space $\R^{2n}/\ell(t)$.

A more explicit description of these curves can be obtained by working with
{\sl frames:}  If $\cA(t)$ is a smooth curve of $2n \times n$ matrices of
rank $n$,  the curve of $n$-dimensional subspaces spanned by the columns of
$\cA(t)$ is fanning if and only if the $2n \times 2n$ matrix
$\join{\cA(t)}{\dot{\cA}(t)}$---formed by juxtaposing $\cA(t)$ and its derivative
$\dot{\cA}(t)$---is invertible for all values of $t$. It will be useful to denote
such curves of $2n \times n$ matrices as {\sl fanning curves of frames,} or as
{\sl fanning frames.} Remark that two fanning frames $\cA(t)$ and $\cB(t)$ span
the same curve of $n$-dimensional subspaces if and only if there is a curve $X(t)$
of invertible $n \times n$ matrices such that $\cB(t) = \cA(t)X(t)$.

\medskip
\noindent {\sl Examples.}
\point If $\cA_1$ and $\cA_2$ are two frames such that the matrix
$\join{\cA_1}{\cA_2}$ is invertible, the {\sl line\/}
$\cA(t) = \cA_1 + t \cA_2$ is a fanning frame.
\point If $\ell$ is an $n$-dimensional subspace of $\R^{2n}$ and $\bX$ is a linear
transformation from $\R^{2n}$ to itself such that $\bX \ell$ is transversal to
$\ell$, the curve $\ell(t) = \exp(t\bX)\ell$ is a fanning curve.
\point P.~Griffiths shows in \cite{Griffiths:1974} that the ruled surface in
real projective space defined by a curve in $\hbox{\rm G}_2$ is developable
if and only if the curve in $\hbox{\rm G}_2$ is {\sl not\/} fanning on any interval.
\medskip

Our main insight is that the key to understanding the geometry of fanning curves is
the following, almost tautological, construction. Given a fanning curve $\ell(t)$,
let $\pi_t$ be the canonical projection from $\R^{2n}$ to $\R^{2n}/\ell(t)$, and
recall that $\dot{\ell}(t)$ is an invertible linear map from $\ell(t)$ to
$\R^{2n}/\ell(t)$. The {\sl fundamental endomorphism\/} $\bF(t)$ at time $t$ is the
linear transformation from  $\R^{2n}$ to itself whose value at a vector $\bv$ is
$(\dot{\ell}(t))^{-1} \pi_t \bv \in \ell(t) \subset \R^{2n}$.
If the curve $\ell(t)$ is spanned by a fanning frame $\cA(t)$, then the matrix for
its fundamental endomorphism in the canonical basis of $\R^{2n}$ is
$$
 \join{\cA(t)}{\dot{\cA}(t)}
\pmatrix{ \mzero & I \cr \mzero & \mzero \cr}
\join{\cA(t)}{\dot{\cA}(t)}^{-1},
$$
where $\mzero$ represents, in this case, the $n \times n$ zero matrix. It follows
immediately from this formula that if $\bT$ is an invertible linear transformation
of $\R^{2n}$ and $\bF(t)$ is the fundamental endomorphism of a fanning curve
$\ell(t)$, the fundamental endomorphism of $\bT\ell(t)$ is $\bT\bF(t)\bT^{-1}$.

Seen as a map that takes one-jets of fanning curves to elements of the Lie algebra
$\lalgebra$ of $2n \times 2n$ matrices, the fundamental endomorphism describes the
prolonged action of $\lgroup{2n}$ on the space of one-jets of fanning curves on the
Grassmannian. Indeed, in Section~7 we shall prove the following characterization.

\newclaim Theorem.
A map from the space of one-jets of fanning curves in $\grassmann$ to the Lie
algebra $\lalgebra$ is equivariant with respect to the $\lgroup{2n}$ action on
these spaces if and only if it is of the form $a\bF + b\bI$, where $\bI$ is the
identity matrix and $a$ and $b$ are real numbers.

The curve $\bF(t)$ associated to a fanning curve $\ell(t)$ has the remarkable
property that its derivative $\dot{\bF}(t)$ is a curve of reflections (i.e.,
$(\dot{\bF}(t))^2 = \bI$) such that the eigenspace of $\dot{\bF}(t)$ associated to
the eigenvalue $-1$ is precisely $\ell(t)$. The eigenspace of $\dot{\bF}(t)$
associated to the eigenvalue $1$ is the {\sl horizontal curve\/} $h(t)$ of $\ell(t)$.

The horizontal curve---appearing for the first time, under a different guise, in
the work of S.~Ahdout \cite{Ahdout:1989}---plays a fundamental role in the
understanding of the geometric invariants of $\ell(t)$. In fact, we will see in
Section~7 that it describes the prolonged action of $\lgroup{2n}$ on the space of
two-jets of fanning curves on the Grassmannian.

\newclaim Theorem.
The assignment that sends a fanning curve $\ell(t)$ to its horizontal curve
$h(t)$ is characterized by the following four properties:
\item{1.} At each time $t$ the subspace $h(t)$ is transversal to $\ell(t)$;
\item{2.} The subspace $h(\tau)$ depends only on the two-jet of the curve $\ell(t)$
          at $t = \tau$.
\item{3.} If $\bT$ is an invertible linear transformation of $\R^{2n}$, the
          horizontal curve of $\bT\ell(t)$ is $\bT h(t)$.
\item{4.} If $\ell(t)$ is spanned by a line $\cA + t\cB$ in the space of frames,
          $h(t)$ is constant.

The main geometric invariant of a fanning curve $\ell(t)$ in the Grassmannian
is its {\sl Jacobi endomorphism\/} $\bK(t) = \ddot{\bF}(t)^2/4$. Alternatively,
if $\bP(t) = (\bI - \dot{\bF}(t))/2$ is the projection onto $\ell(t)$ with kernel
$h(t)$, $\bK(t) = \dot{\bP}(t)^2$. Hence, the Jacobi endomorphism describes how
the horizontal curve moves with respect to $\ell(t)$.

\newclaim Theorem.
Let $\ell(t)$ be a fanning curve in the Grassmannian $\grassmann$ and let $h(t)$ be
its horizontal curve. The Jacobi endomorphism $\bK(t)$ of $\ell(t)$ satisfies the
following properties:
\item{1.} At each time $t$ the endomorphism $\bK(t)$ preserves the decomposition
          $\R^{2n} = \ell(t) \oplus h(t)$.
\item{2.} If $\bT$ is an invertible linear map from $\R^{2n}$ to itself, the Jacobi
          endomorphism of $\bT \ell(t)$ is $\bT\bK(t)\bT^{-1}$.
\item{3.} If $s$ is a diffeomorphism of the real line and $\{s(t),t\}$ denotes its
Schwarzian derivative, the Jacobi endomorphism of $\ell(s(t))$ is
$$
\bK(s(t))\dot{s}(t)^2 + (1/2)\{s(t),t\}\,\bI .
$$

The Jacobi endomorphism clarifies the geometry underlying the Schwarzian derivative
and its matrix generalizations. If $\cA(t)$ is a fanning frame, at each instant
$t$ the columns of $\cA(t)$ and $\dot{\cA}(t)$ span $\R^{2n}$ and, therefore, we
have the differential equation
$$
\ddot{\cA} + \dot{\cA}P(t) + \cA Q(t) = \mzero, \leqno{(1.1)}
$$
where $Q(t)$ and $P(t)$ are smooth curves of $n \times n$ matrices. Let us
define the {\sl Schwarzian\/} of $\cA(t)$ as the function
$$
\{\cA(t),t\} = 2Q(t) - (1/2)P(t)^2 - \dot{P}(t). \leqno{(1.2)}
$$
When the fanning frame is of the form
$$
\cA(t) = \pmatrix{I \cr M(t) \cr}, \quad \hbox{\rm then} \quad
\{\cA(t),t\} = \frac{d}{dt}(\dot{M}^{-1}\ddot{M}) -
(1/2)(\dot{M}^{-1}\ddot{M})^2,
$$
is the matrix Schwarzian introduced by B.~Schwarz~\cite{Schwarz:1979} and
M.~I.~Zelikin~\cite{Zelikin:1992}.

\newclaim Theorem.
The Schwarzian of a fanning frame $\cA(t)$ is characterized by the
equation $\bK(t)\cA(t) = (1/2)\cA(t)\{\cA(t),t\}$.

The properties of the Schwarzian follow immediately from those of the Jacobi
endomorphism.

\newclaim Corollary.
The Schwarzian of a fanning frame $\cA(t)$ satisfies the following
properties:
\item{1.} If $\bT$ is an invertible linear transformation from $\R^{2n}$ to itself,
$\{\bT \cA(t),t\} = \{\cA(t),t\}$.
\item{2.} If $X(t)$ is a curve of invertible $n \times n$ matrices,
$\{\cA(t)X(t),t\} = X(t)^{-1}\{\cA(t),t\}X(t)$.
\item{3.} If $s$ is a diffeomorphism of the real line,
$\{\cA(s(t)),t\} = \{A(s(t)),s\}\dot{s}(t)^2 + \{s(t),t\}I$.

Because of Property~(2), the only thing we can say about fanning frames
spanning congruent curves in $\grassmann$ is that their
Schwarzians are pointwise conjugate. In order to obtain a deeper understanding
of the geometry of fanning curves, in Section~4 we introduce {\sl normal frames}
for which Equation~(1.1) takes the normal form
$$
\ddot{\cA} + (1/2) \cA\, \{\cA(t),t\} = \mzero  . \leqno{(1.3)}
$$

Using normal frames, the ``fundamental theorem'' for fanning curves on the
Grassmannian is an easy consequence of the uniqueness theorem for solutions of
differential equations.

\newclaim Theorem.
Two fanning curves  of $n$-dimensional subspaces of $\R^{2n}$ are congruent if
and only if the Schwarzians of any two of their normal frames are conjugate by
a constant $n \times n$ invertible matrix.

The fundamental theorem for curves of Lagrangian subspaces under the action of
the linear symplectic group is only slightly more involved and will be presented
in Section~6.

The Jacobi endomorphism is also central to the study of unparameterized fanning
curves in $\grassmann$. Again, this is a classical subject---when $n=2$ this is
the projective differential geometry of nondevelopable ruled surfaces in
$\RP^3$ (see, for example, \cite{Wilczynski:1906})---, but the new approach
eliminates many of the computations and gives a better understanding of the subject.

It is clear from Property~(3) in Theorem~1.4 that any fanning curve admits a
{\sl special parameterization\/} so that the trace of its Jacobi endomorphism
vanishes identically. Moreover, any two such parameterizations are projectively
equivalent. It follows that in a special parameterization the operator-valued
quadratic differential $\bK(t)dt^2$ is an invariant of the unparameterized curve.
Defining a {\sl special normal frame\/} as a normal frame that spans a
fanning curve with a special parameterization, the fundamental theorem for
unparameterized fanning curves may be stated as follows:

\newclaim Theorem.
Two unparameterized fanning curves  of $n$-dimensional subspaces of $\R^{2n}$ are
congruent if and only if up to a projective change of parameters the Schwarzians of
any two of their special normal frames are conjugate by a constant $n \times n$
invertible matrix.

In the case $n=2$ this result is due to Wilczynski
(see \cite[pp. 114--116]{Wilczynski:1906}). Indeed, in some sense the
fundamental theorems for fanning curves in $\grassmann$ and $\lagrangian$,
parameterized or not, are just geometric reinterpretations of old work
on the invariant theory of linear ordinary differential equations. We have
included them in this paper because after a long search of the literature we were
unable to find suitable references.

Also in this paper the reader will find a characterization of fanning curves in
$\grassmann$ that are projections of one-parameter subgroups of the linear group
$\lgroup{2n}$ (Section~5), a thorough study of the geometry of fanning curves
of Lagrangian subspaces (Section~6), and a simplified account of the approaches
of S.~Ahdout and A.~A.~Agrachev et al. to the geometry of fanning curves (Section~8).
In our opinion, the additional insight that a comparisons of approaches will give
the reader warrants this small scholarly effort on our part. We also think this will
help break the cycle of rediscovery in which the subject has been caught up for
some time.

The close relationship between the geometry of fanning curves in the Grassmannian
and the elegant formalisms of connections for second-order differential equations
(semi-sprays) developed by J.~Klein and A.~Voutier \cite{Klein-Voutier:1968},
J.~Grifone~\cite{Grifone:1972a}, and P.~Foulon~\cite{Foulon:1986} will be the
subject of a future publication. At the end of this lengthy introduction we
can succinctly summarize the contents of this paper and of that future publication
as follows: {\sl the geometric invariants of fanning curves in the Grassmannian
on $n$-dimensional subspaces of $\R^{2n}$ arise from the fundamental
endomorphism $\bF$---a tautological construct---and its derivatives $\dot{\bF}$
and $\ddot{\bF}$, from which we define the Jacobi endomorphism
$\bK = (1/4)\ddot{\bF}{}^2$. Using the construction outlined in the first paragraph
of the introduction, these correspond in a precise way to the vertical endomorphism,
connection, and curvature (Jacobi endomorphism) introduced by Klein, Voutier,
Grifone and Foulon in their study of connections for semi-sprays and Finsler
metrics.}

\mysection{The fundamental endomorphism and its derivatives}
The two groups acting naturally on the space of fanning curves in $\grassmann$ are
the group of invertible linear transformations of $\R^{2n}$ and the group
of diffeomorphisms of the real line acting by reparameterizations. If we
wish to work with fanning frames, we must add the group  of
smooth curves of $n \times n$ invertible matrices acting by
$(\cA(t),X(t)) \mapsto \cA(t)X(t)$. Any quantity associated to a fanning frame
$\cA(t)$ that is invariant under this action depends only on the fanning curve
on the Grassmannian defined by the span of the columns of $\cA(t)$.

We caution the reader that unless we specifically state otherwise the canonical
basis of $\R^{2n}$ shall be freely used to identify linear transformations of
$\R^{2n}$ with $2n \times 2n$ matrices. Despite this, all constructions are
intrinsic in nature and apply to curves of $n$-dimensional linear subspaces in a
$2n$-dimensional vector space over the real numbers.

\newpar Definition.
The {\sl fundamental endomorphism\/} of a fanning frame $\cA(t)$ at
a given time $\tau$ is the linear transformation from $\R^{2n}$ to itself
defined by the equations $\bF(\tau)\cA(\tau) = \mzero$,
$\bF(\tau)\dot{\cA}(\tau) = \cA(\tau)$.

Equivalently, as in the introduction, we could have defined the fundamental
endomorphism by the formula
$$
\bF(t) = \join{\cA(t)}{\dot{\cA}(t)}
\pmatrix{ \mzero & I \cr \mzero & \mzero \cr}
\join{\cA(t)}{\dot{\cA}(t)}^{-1}.
$$
Using this formula, an easy calculation suffices to establish the main properties
of the fundamental endomorphism.

\newclaim Proposition.
Let $\cA(t)$ be a fanning frame. Its fundamental endomorphism $\bF(t)$
satisfies the following properties:
\item{1.} If $X(t)$ is a smooth curve of invertible $n \times n$ matrices, the
fundamental endomorphism of $\cA(t)X(t)$ is again $\bF(t)$.
\item{2.} If $\bT$ is an invertible linear transformation from $\R^{2n}$ to itself,
the fundamental endomorphism of $\bT \cA(t)$ is $\bT\bF(t)\bT^{-1}$.
\item{3.} The fundamental endomorphism of the reparameterized curve
          $\cA(s(t))$ is $\bF(s(t))\dot{s}(t)^{-1}$.

By the first property, the fundamental endomorphism is defined for fanning
curves in the Grassmannian. The intrinsic definition given in the
introduction---suggested to us by Fran Burstall---justifies the term
``fundamental''. On the other hand, basic properties such as the one given in
the following proposition are hard to prove intrinsically.

\newclaim Proposition.
Let  $\bF(t)$ be the fundamental endomorphism of a fanning frame
$\cA(t)$.  At each value of $t$ the derivative $\dot{\bF}(t)$ is a reflection
$(\hbox{i.e.}, \ (\dot{\bF}(t))^2 = \bI))$ whose $-1$ eigenspace is spanned by the
columns of $\cA(t)$.

\Proof
We first show that $\dot{\bF}(t) \cA(t) = -\cA(t)$. From this follows that
$\dot{\bF}(t)$ restricted to the subspace $\ell(t)$ spanned by the columns of
$\cA(t)$ is minus the identity. Since $\ell(t)$ is the range of
$\bF(t)$, this also implies that $\dot{\bF}(t)\bF(t) = - \bF(t)$.
Differentiating the equation $\bF(t)\cA(t) = \mzero$ and recalling that
$\bF(t)\dot{\cA}(t) =  \cA(t)$,  it follows that
$\dot{\bF}(t)\cA(t) + \bF(t) \dot{\cA}(t) = \mzero$ and
$\dot{\bF}(t)\cA(t) = -\cA(t)$.

We now show that $(\dot{\bF}(t))^{2}\dot{\cA}(t) = \dot{\cA}(t)$. Since we already
know that $(\dot{\bF}(t))^{2}\cA(t) = \cA(t)$, this will prove that $\dot{\bF}(t)$
is a reflection for every value of $t$.
Differentiating the equation $\bF(t)\dot{\cA}(t) = \cA(t)$ and using that
$\dot{\bF}(t)\bF(t) = - \bF(t)$, we have that
$$
(\dot{\bF}(t))^2\dot{\cA}(t) = \dot{\bF}(t)(\dot{A}(t) - \bF(t)\ddot{\cA}(t))
=\dot{A}(t) . \qed
$$

We remark that since the fundamental endomorphism depends only on the curve on
the Grassmannian and not on the fanning frame used to represent it,
the same holds for all its time derivatives. In particular, the curve of
reflections $\dot{\bF}(t)$ and the associated curve of projections
$\bP(t) = (\bI - \dot{\bF}(t))/2$ depend only on the curve on the Grassmannian.

\newpar Definition.
Let $\ell(t)$ be a fanning curve in the Grassmannian and let $\bF(t)$ be its
fundamental endomorphism. The map that takes $t$ to the  kernel of the projection
$\bP(t) = (\bI - \dot{\bF}(t))/2$ is the {\sl horizontal curve\/} of
$\ell(t)$.

It is clear from the definition that the subspace $h(\tau)$ is transversal to
$\ell(\tau)$ and depends only on the two-jet of $\ell(t)$ at $t = \tau$.

\newclaim Proposition.
Let $\ell(t)$ be a fanning curve on $\grassmann$ and let $h(t)$ be its horizontal curve.
If $\bT$ is an invertible linear transformation from $\R^{2n}$ to itself, the
horizontal curve of $\bT \ell(t)$ is $\bT h(t)$.

\Proof
This follows immediately from the fact that if $\bT$ is an invertible linear
transformation from $\R^{2n}$ to itself, the fundamental endomorphism of
$\bT \ell(t)$ and its derivative are $\bT\bF(t)\bT^{-1}$ and
$\bT\dot{\bF}(t)\bT^{-1}$, respectively. \qed

We now turn to the study of the second derivative $\ddot{\bF}$ of the
fundamental endomorphism. The geometric meaning of the computations
will be clearer if we work with $\dot{\bP} = -(1/2)\ddot{\bF}$.

\newclaim Proposition.
Let $\ell(t)$ be a fanning curve in the Grassmannian $\grassmann$ and let $h(t)$ be
its horizontal curve. If $\bP(t)$ denotes the projection onto $\ell(t)$ with
kernel $h(t)$, then $\dot{\bP}(t)$ maps $h(t)$ into $\ell(t)$ and maps $\ell(t)$
isomorphically to $h(t)$.

\Proof
Differentiating the identity $\bP(t)^2 = \bP(t)$ we have that
$\dot{\bP}(t)\bP(t) = (\bI - \bP(t))\dot{\bP}(t)$. If we notice that
$\bI - \bP(t)$ is the projection onto $h(t)$ with kernel $\ell(t)$,
the equation
$$
\dot{\bP}(t)(\ell(t)) = \dot{\bP}(t)\bP(t)(\ell(t)) =
(\bI - \bP(t))\dot{\bP}(t)(\ell(t))
$$
implies that the subspace $\dot{\bP}(t)(\ell(t))$ is contained in $h(t)$.
The proof that the subspace $\dot{\bP}(t)(h(t))$ is contained in $\ell(t)$
is nearly identical.

In order to verify that $\dot{\bP}(t)$ maps $\ell(t)$ isomorphically onto
$h(t)$, we make use of the identity $\dot{\bP} = (-1/2)\ddot{\bF}$.
It follows from the proof of Proposition~2.3 that if $\cA(t)$ is a fanning
frame spanning $\ell(t)$, then $\dot{\bF}(t) \cA(t) = -\cA(t)$ and
$\dot{\bF}(t) \dot{\cA}(t) = \dot{\cA}(t) - \bF(t)\ddot{\cA}(t)$.
Differentiating the first of these equations and using the second, we obtain
$\dot{\bP}(t) \cA(t) = \dot{\cA}(t) - (1/2)\bF(t)\ddot{\cA}(t)$.
Since the columns of $\bF(t)\ddot\cA(t)$ are linear combinations of those of
$\cA(t)$, and $\join{\cA(t)}{\dot{\cA}(t)}$ is invertible, it follows that
$\dot{\bP}(t) \cA(t)$ has rank $n$ and must span $h(t)$. \qed

A remark on the preceding proof is that $\dot{\bP}(t) \cA(t)$
is the projection of $\dot{\cA}(t)$ onto the horizontal subspace $h(t)$.
Indeed, on differentiating the equality $\bP(t)\cA(t) = \cA(t)$, we obtain
$\dot{\bP}(t) \cA(t) = (\bI - \bP(t))\dot{\cA}(t)$.

\newpar Definition.
The {\sl horizontal derivative\/} of a  fanning frame $\cA(t)$
is the curve of frames
$$
\cH(t) := (\bI - \bP(t))\dot{\cA}(t) = \dot{\bP}(t)\cA(t) =
\dot{\cA}(t) - (1/2)\bF(t)\ddot{\cA}(t) = -(1/2)\ddot{\bF}(t)\cA(t).
$$

\noindent {\sl Remark.}
If, as in Equation~(1.1), we write $\ddot{\cA} + \dot{\cA}P(t) + \cA Q(t) = \mzero$,
then the horizontal derivative of $\cA(t)$ is
$$
\cH(t) = \dot{\cA}(t) + (1/2)\cA(t)P(t). \leqno{(2.1)}
$$

From the horizontal derivative we immediately obtain an elementary description
of the horizontal subspace of a fanning curve.

\newclaim Proposition.
Let $\ell(t)$ be a fanning curve on the Grassmannian. If $\cA_{\tau}(t)$
is a fanning frame spanning $\ell(t)$ that satisfies
$\ddot{\cA}_\tau(\tau) = \mzero$, then the columns of $\dot{\cA}_\tau(\tau)$
span the horizontal subspace of $\ell(t)$ at $t = \tau$.

In the two-dimensional case this implies that if $\ell(t)$ is a fanning curve
in the projective line and $h(t)$ is its horizontal curve, then a straight line
$l_\tau$ on the plane not passing through the origin is parallel to $h(\tau)$ if
and only if the acceleration of the curve of vectors
$\bv_{\tau}(t) = \ell(t) \cap l_\tau$ is zero at $t = \tau$.

An easy calculation shows that the horizontal derivative is well-behaved with
respect to the three natural group actions on the space of fanning frames.

\newclaim Proposition.
The horizontal derivative $\cH(t)$ of a fanning frame $\cA(t)$ satisfies
the following properties:
\item{1.} If $X(t)$ is a smooth curve of $n \times n$ invertible matrices,
the horizontal derivative of $\cA(t)X(t)$ is $\cH(t)X(t)$.
\item{2.} If $\bT$ is an invertible linear transformation from $\R^{2n}$ to
itself, the horizontal derivative of $\bT \cA(t)$ is $\bT \cH(t)$.
\item{3.} The horizontal derivative of the reparameterized curve $\cA(s(t))$ is
$$
\cH(s(t))\dot{s}(t) + (1/2)\cA(s(t))\dot{s}(t)^{-1}\ddot{s}(t) .
$$

The horizontal derivative is useful in many computations like the one in the
proof of the following interesting property of the fundamental endomorphism
and its derivatives.

\newclaim Proposition.
Let $\ell(t)$ be a fanning curve with fundamental endomorphism $\bF(t)$. If
$[\bA,\bB]_+$ denotes the anticommutator $\bA\bB + \bB\bA$, then
$[\bF(t),\dot{\bF}(t)]_+ = \mzero$, $[\dot{\bF}(t),\ddot{\bF}(t)]_+ = \mzero$,
and $[\bF(t),\ddot{\bF}(t)]_+ = -2\bI$.

\Proof
The first two identities are obtained by differentiating the identities
$\bF(t)^2 = \mzero$ and $\dot{\bF}(t)^2 = \bI$. To obtain the third we
take a curve of frames $\cA(t)$ spanning $\ell(t)$ and
compute both $[\bF(t),\ddot{\bF}(t)]_+\, \cA(t)$ and
$[\bF(t),\ddot{\bF}(t)]_+\, \cH(t)$.
In the first of these anticommutators, notice that $\ddot{\bF}\bF\cA = \mzero$
and that
$$
\bF\ddot{\bF}\cA = \bF(-2\cH) = \bF(-2\dot{\cA} + \bF\ddot{\cA}) = -2\cA .
$$
In the second anticommutator, we have that
$\bF\ddot{\bF}\cH = -2\bF\dot{\bP}\cH = \mzero$
because $\dot{\bP}$ sends $h(t)$ into $\ell(t)$, the kernel of $\bF(t)$.
Moreover,
$$
\ddot{\bF}\bF\cH = \ddot{\bF}\bF(\dot{\cA} - (1/2)\bF\ddot{\cA}) =
\ddot{\bF}{\cA} = -2\cH .
$$
We have then that $[\dot{\bF}(t),\ddot{\bF}(t)]_+\, \cA(t) = -2\cA(t)$ and
that $[\dot{\bF}(t),\ddot{\bF}(t)]_+\, \cH(t) = -2\cH(t)$, which means that
$[\dot{\bF}(t),\ddot{\bF}(t)]_+ = -2\bI$. \qed

\mysection{The Jacobi endomorphism and the Schwarzian}
The natural notion of curvature for fanning curves in the Grassmannian $\grassmann$
is given by an operator-valued function of the parameter that is closely related
to the Schwarzian derivative and its matrix generalizations.

\newpar Definition.
Let $\ell(t)$ be a fanning curve in the Grassmannian $\grassmann$, let $\bF(t)$ be
its fundamental endomorphism, and let $h(t)$ be its horizontal curve. The
{\sl Jacobi endomorphism\/} of $\ell(t)$ is $\bK(t) = \ddot{\bF}(t)^2/4$.
Alternatively, if $\bP(t) = (\bI - \dot{\bF}(t))/2$ is the projection onto
$\ell(t)$ with kernel $h(t)$, $\bK(t) = \dot{\bP}(t)^2$.

The intuition behind the Jacobi endomorphism is that measures how the horizontal
curve moves with respect to $\ell(t)$. More precisely, if $\cA(t)$ is a fanning
frame spanning $\ell(t)$ and $\cH(t)$ is its horizontal derivative,
$$
\bP(t)\dot{\cH}(t) = -\dot{\bP}(t)\cH(t) = -\dot{\bP}(t)^2 \cA(t)
= - \bK(t) \cA(t) . \leqno{(3.1)}
$$

\newclaim Theorem.
Let $\ell(t)$ be a fanning curve in $\grassmann$ and let $h(t)$ be its horizontal
curve. The Jacobi endomorphism $\bK(t)$ of $\ell(t)$ satisfies the following
properties:
\item{1.} At each time $t$ the endomorphism $\bK(t)$ preserves the decomposition
          $\R^{2n} = \ell(t) \oplus h(t)$.
\item{2.} If $\bT$ is an invertible linear map from $\R^{2n}$ to itself, the Jacobi
          endomorphism of $\bT \ell(t)$ is $\bT\bK(t)\bT^{-1}$.
\item{3.} If $s$ is a diffeomorphism of the real line, the Jacobi endomorphism of
$\ell(s(t))$ is
$$
\bK(s(t))\dot{s}(t)^2 + (1/2)\{s(t),t\}\,\bI .
$$

\Proof
The first property follows immediately from the identity $\bK(t) = \dot{\bP}(t)^2$
and Proposition~2.6. The second property follows from the $\lgroup{2n}$-equivariance
of the fundamental endomorphism.

By Proposition~2.2, if $\bF(t)$ denotes the fundamental endomorphism of $\ell(t)$,
the fundamental endomorphism of $\ell(s(t))$ is $\bF(s(t))\dot{s}(t)^{-1}$.
Setting, for convenience, $r(t) = \dot{s}(t)^{-1}$, we have that the second
derivative of the fundamental endomorphism of $\ell(s(t))$ is
$$
\bF''(s(t))\dot{s}(t) + \bF'(t)\dot{s}(t)\dot{r}(t) + \bF(s(t))\ddot{r}(t),
\quad \hbox{\rm where} \quad {}' = d/ds .
$$

Writing the square of this expression in terms of anticommutators and applying
the identities of Proposition~2.10 together with the identities $\bF^2 = \mzero$
and $(\bF')^2 = \bI$, we obtain
$$
\left(\frac{d^2}{dt^2}\, \bF(s(t))\dot{s}(t)^{-1} \right)^2 =
\bF''(s(t))^2 \dot{s}(t)^2 +
\bI (\dot{s}(t)^2 \dot{r}(t)^2 - 2\dot{s}(t)\ddot{r}(t)).
$$
A short calculation reveals that
$(\dot{s}(t)^2 \dot{r}(t)^2 - 2\dot{s}(t)\ddot{r}(t))$ is twice
the Schwarzian $\{s(t),t\}$. We conclude that the Jacobi endomorphism of
$\ell(s(t))$ is $\bK(s(t))\dot{s}(t)^2 +  (1/2)\{s(t),t\}\,\bI$ . \qed

\newpar Definition.
If $\cA(t)$ is a fanning frame and the curves of $n \times n$
matrices $Q(t)$ and $P(t)$ are defined by the equation
$\ddot{\cA} + \dot{\cA}P(t) + \cA Q(t) = \mzero$, the function
$$
\{\cA(t),t\} = 2Q(t) - (1/2)P(t)^2 - \dot{P}(t)
$$
is the {\sl Schwarzian\/} of $\cA(t)$.

\newclaim Theorem.
The Schwarzian of a fanning frame $\cA(t)$ is characterized by the
equation $\bK(t)\cA(t) = (1/2)\cA(t)\{\cA(t),t\}$.

\Proof
The proof follows from the equations $\bK(t)\cA(t) = -\bP(t)\dot{\cH}(t)$ and
$\cH(t) = \dot{\cA}(t) + (1/2)\cA(t)P(t)$. Indeed, on differentiating this
last equation and replacing $\ddot{\cA}(t)$ by $-\dot{\cA}(t)P(t) -\cA(t)Q(t)$,
we obtain
$$
\eqalign{
\dot{\cH}(t) &= \cA(t)(-Q(t) + (1/4)P(t)^2 + (1/2)\dot{P}(t)) - (1/2)\cH(t)P(t) \cr
             &= -(1/2)\cA(t)\{\cA(t),t\} - (1/2)\cH(t)P(t). \cr}
\leqno{(3.2)}
$$
Applying $-\bP(t)$ to both sides of the equation yields the result. \qed

\newclaim Corollary.
If $\cA(t)$ is a fanning frame and $\cH(t)$ is its horizontal derivative,
the matrix for the Jacobi endomorphism $\bK(t)$ of $\cA(t)$ in the basis of
$\R^{2n}$ formed by the columns of $\join{\cA(t)}{\cH(t)}$ is
$$
\pmatrix{ \{\cA(t),t\}/2 & \mzero \cr \mzero & \{\cA(t),t\}/2 \cr} .
$$

\Proof
Since we already know that $\bK(t)\cA(t) = (1/2)\cA(t)\{\cA(t),t\}$, it remains
to compute $\bK(t)\cH(t) = \dot{\bP}(t)^2\cH(t)$.
By Equation~(3.1),
$\dot{\bP}\cH(t) = \bK(t)\cA(t) = (1/2)\cA(t)\{\cA(t),t\}$. This yields
$$
\bK(t)\cH(t) = \dot{\bP}(t)^2\cH(t) = (1/2)\dot{\bP}(t) \cA(t)\{\cA(t),t\}
= (1/2)\cH(t)\{\cA(t),t\}.
$$
\qed

The properties of the Schwarzian follow immediately from those  of the Jacobi
endomorphism.

\newclaim Corollary.
The Schwarzian of a fanning frame $\cA(t)$ satisfies the following
properties:
\item{1.} If $\bT$ is an invertible linear transformation from $\R^{2n}$ to itself,
$\{\bT \cA(t),t\} = \{\cA(t),t\}$.
\item{2.} If $X(t)$ is a curve of invertible $n \times n$ matrices,
$\{\cA(t)X(t),t\} = X(t)^{-1}\{\cA(t),t\}X(t)$.
\item{3.} If $s$ is a diffeomorphism of the real line,
$\{\cA(s(t)),t\} = \{A(s(t)),s\}\dot{s}(t)^2 + \{s(t),t\}I$.
\item{4.} If the fanning frame is of the form
$$
\cA(t) = \pmatrix{I \cr M(t) \cr}, \quad \hbox{\rm then} \quad
\{\cA(t),t\} = \frac{d}{dt}(\dot{M}^{-1}\ddot{M}) -
(1/2)(\dot{M}^{-1}\ddot{M})^2
$$

As a simple application of these properties, we obtain the precise---and apparently
new---transformation law for the matrix Schwarzian under matrix linear fractional
transformations (compare with Proposition~6.14 in \cite[p. 207]{Zelikin:2000}).

\newclaim Corollary.
Let $A$, $B$, $C$, and $D$ be $n \times n$ matrices, and let
$M(t)$ be a smooth curve of $n \times n$ matrices such that its derivative
$\dot{M}(t)$ is never singular. If $(A + BM(t))$ is invertible
for all values of $t$ in some interval and if $S_t(M)$ denotes the matrix
Schwarzian of $M(t)$, then the matrix Schwarzian of $(C + DM(t))(A + BM(t))^{-1}$
on that interval is $(A + BM(t))S_t(M) (A + BM(t))^{-1}$.

\Proof
By the preceding corollary, the matrix Schwarzian of $(C + DM(t))(A + BM(t))^{-1}$
is the Schwarzian of the fanning frame
$$
\cA(t) = \pmatrix{A + BM(t) \cr C + DM(t) \cr} (A + BM(t))^{-1}
= \pmatrix{A & B \cr C & D \cr} \pmatrix{I \cr M(t) \cr} (A + BM(t))^{-1} .
$$
The result now follows from the first two properties of the Schwarzian of
frames (Corollary~3.6). \qed

\subsection{Reparameterization invariants of  fanning curves}
The reader has undoubtedly noticed that the traces of powers and
derivatives of the Jacobi endomorphism of a fanning curve do not
change if the curve is subjected to the action of linear transformations
on the Grassmannian $\grassmann$. However, these quantities are not
invariant under reparameterizations of the curve.

The simplest quantity associated to a fanning curve that is invariant
under both reparameterizations and the action of linear transformations on the
Grassmannian $\grassmann$ is the fourth-order differential
$$
\left(\tr \, \bK(t)^2 - (1/n)(\tr\, \bK(t))^2\right)(dt)^4 .
$$
Under a different notation, this invariant appeared for the first time
in Agrachev and Zelenko~\cite{Agrachev-Zelenko:2002a}

Another interesting invariant (or rather covariant) is given by the
operator-valued, fifth-order differential $[\dot{\bK}(t),\bK(t)](dt)^5$ which is
the analogue for fanning curves of a projective invariant for semi-sprays
that is described by Foulon in Proposition~VI.5 of \cite{Foulon:1986}.

The proof of the invariance of the fourth and fifth order differentials
just defined is a simple computation involving Property~(3) in Theorem~3.2.
However, the most important information that can be extracted from this property
is that any fanning curve admits a {\sl special parameterization.}

\newpar Definition.
A fanning curve $\ell(t)$ is said to be {\sl specially parameterized\/}
if the trace of its Jacobi endomorphism vanishes identically.

Two special parameterizations are related by a projective transformation and,
therefore, in a special parameterization the operator-valued quadratic
differential $\bK(t)dt^2$ is an invariant of the unparameterized curve. We shall
see in the next section in what way this is the "fundamental invariant" for
unparameterized fanning curves in $\grassmann$.

\mysection{Normal frames and solution of the congruence problem}
As we have remarked all along, if a fanning curve $\ell(t)$ on the Grassmannian
is spanned by a frame $\cA(t)$ it is also spanned by all frames of the form
$\cA(t)X(t)$, where $X(t)$ is a smooth curve of invertible $n \times n$ matrices.
Since $\{\cA(t)X(t),t\} = X(t)^{-1} \{\cA(t),t\}X(t)$, the Schwarzians of any
two frames for $\ell(t)$ are point-wise conjugate. Hence, quantities such as the
traces of the Schwarzian and its powers will depend only on the fanning curve on
the Grassmannian and are, furthermore, invariant under the action of the linear
group on the space of fanning curves. Nevertheless, if we wish to obtain a deeper
understanding of the geometry of a fanning curve, we need to introduce a class of
frames that is better adapted to it.

\newpar Definition.
A fanning frame $\cA(t)$ is said to be {\sl normal\/} if the columns
of its second derivative $\ddot{\cA}(t)$ are linear combinations of the columns
of $\cA(t)$ for all values of $t$.

\newclaim Proposition.
If $\ell(t)$ is a fanning curve of $n$-dimensional subspaces in $\R^{2n}$, there
exists a normal frame that spans it. Moreover, if $\cA(t)$ and $\cB(t)$ are
two normal frames spanning $\ell(t)$, there exists a fixed invertible $n \times n$
matrix $X$ such that $\cB(t) = \cA(t)X$.

\Proof
Let $\cA(t)$ be any fanning frame spanning $\ell(t)$ and let $P(t)$ be the
curve of $n \times n$ matrices defined by the equation
$\ddot{\cA}(t) + \dot{\cA}(t)P(t) + \cA(t)Q(t) = \mzero$. If
$X(t)$ is the the curve of $n \times n$ matrices that solves the initial value
problem $\dot{X}(t) = (1/2)P(t)X(t)$, $X(0) = I$, an
easy computation shows that $\cA(t)X(t)$ is a normal frame. Notice that since $X(0)$
is invertible, $X(t)$ is invertible for all values of $t$ and $\cA(t)X(t)$
is again a fanning frame.

If $\cA(t)$ and $\cB(t)$ are fanning frames spanning the same curve $\ell(t)$
in $\grassmann$, then there exists a curve of $n \times n$ invertible matrices
$X(t)$ such that $\cB(t) = \cA(t)X(t)$. Differentiating this equation twice and
using that $\cA(t)$ is normal, we see that the only possible way in which the
columns of $\ddot{\cB}(t)$ could be linear combinations of the columns of
$\cB(t)$ is that $\dot{X}(t)$ be identically zero. \qed

The use of normal frames leads to a simple solution of the congruence problem
for (parameterized) fanning curves:
given two fanning curves $\ell(t)$ and $\tilde{\ell}(t)$ in $\grassmann$, when
does there exist an invertible linear transformation $\bT$ of $\R^{2n}$ such
that $\tilde{\ell}(t) = \bT \ell(t)$?

\newclaim Theorem.
Two fanning curves  of $n$-dimensional subspaces of $\R^{2n}$ are congruent if
and only if the Schwarzians of any two of their normal frames are conjugate by
a constant $n \times n$ invertible matrix.

\Proof
If $\cA(t)$ and $\cB(t)$ are two normal frames spanning congruent fanning curves,
then there exists a linear transformation $\bT$ of $\R^{2n}$ such that $\bT \cA(t)$
and $\cB(t)$ span the same curve. Since $\bT\cA(t)$ is again a normal frame,
Corollary~3.6 together with Proposition~4.2 tell us that the Schwarzians of
$\cA(t)$ and $\cB(t)$ are conjugate by a constant $n \times n$ invertible matrix.

On the other hand, if $\cA(t)$ and $\cB(t)$ are two normal frames such that
$X^{-1}\{\cA(t),t\}X = \{\cB(t),t\}$, the normal frame $\cA(t)X$  spans
the same curve as $\cA(t)$ and has the same Schwarzian as $\cB(t)$. Without
loss of generality we could have then assumed that $\cA(t)$ and $\cB(t)$ are
two normal frames with the same Schwarzian.

In order to show that these frames are congruent, we set
$\bT = \join{\cB(0)}{\dot{\cB}(0)} \join{\cA(0)}{\dot{\cA}(0)}^{-1}$
and remark that $\cD(t) = \bT\cA(t)$ satisfies the same
second-order differential equation, $\ddot{\cD} + (1/2)\cD\{\cB(t),t\} = \mzero$, as
$\cB(t)$ and has the same initial conditions. It follows that $\cD(t) = \cB(t)$
and, therefore, $\cA(t)$ is congruent to $\cB(t)$. \qed

\subsection{Fundamental theorem for unparameterized fanning curves}
In order to state and prove the analogue of Theorem~4.3 for unparameterized
fanning curves in $\grassmann$, we define a {\sl special normal frame\/} as a
a normal frame whose Schwarzian has zero trace.

\newclaim Proposition.
If $\ell(t)$ is a fanning curve of $n$-dimensional subspaces in $\R^{2n}$, there
exists a special normal frame that spans a reparameterization of it. Moreover, if
$\cA(t)$ and $\cB(t)$ are special normal frames spanning the same unparameterized
fanning curve in $\grassmann$, there exist a fixed invertible $n \times n$ matrix
$X$ and a projective transformation $\sigma$ of the real line such that
$\cB(\sigma(t)) = \cA(t)X(t)$.

\Proof
Let $\ell(t(s))$ be a special parameterization of the fanning curve traced by
$\ell(t(s))$. The trace of the Jacobi endomorphism of $\ell(t(s))$ is then
identically zero and, by Corollary~3.5, so is the trace of the Schwarzian of
any frame spanning it. It follows than any normal frame spanning $\ell(t(s))$ is
a special normal frame. The question of uniqueness is addressed in a way with
which by now the reader is familiar.
\qed

\newclaim Theorem.
Two unparameterized fanning curves  of $n$-dimensional subspaces of $\R^{2n}$ are
congruent if and only if up to a projective change of parameters the Schwarzians of
any two of their special normal frames are conjugate by a constant $n \times n$
invertible matrix.

The proof is almost identical to that of Theorem~4.3.

\mysection{Special classes of fanning curves}
As a second application of the concept of normal frames we prove the following
characterization of fanning curves with zero Jacobi endomorphism (compare with
Theorem~1 in \cite{Agrachev-Zelenko:2002b}).

\newclaim Theorem.
The following conditions on a fanning curve $\ell(t)$ of $n$-dimensional subspaces
in $\R^{2n}$ are equivalent.
\item{1.} The Jacobi endomorphism of $\ell(t)$ is zero.
\item{2.} The Schwarzian of any curve of frames spanning $\ell(t)$ is zero.
\item{3.} Any normal frame spanning $\ell(t)$ is a line in the space of frames.
\item{4.} The horizontal curve of $\ell(t)$ is constant.

\Proof
The equivalence of (1) and (2) follows immediately from Theorem~3.4. That (2)
implies (3) follows from Equation~(1.3). Indeed, if $\cA(t)$ is a normal
frame $\ddot{\cA}(t) = -(1/2)\cA(t)\{\cA(t),t\}$. If the Schwarzian is zero,
then $\cA(t) = \cA(0) + t \dot{\cA}(0)$.

Since the horizontal derivative of a line $\cA(t) = \cA(0) + t \dot{\cA}(0)$
is constantly equal to $\dot{\cA}(0)$, (3) clearly implies (4). In order to see
that (4) implies (2), thereby finishing the proof, notice that if the horizontal
curve of $\ell(t)$ is constant, the horizontal derivative $\cH(t)$ of any fanning
frame $\cA(t)$ spanning $\ell(t)$ is of the form $\cH(t) = \cH(0)X(t)$ for some
curve $X(t)$ of $n \times n$ invertible matrices. In this case
$\dot{\cH}(t) = \cH(0)\dot{X}(t) = \cH(t)X(t)^{-1}\dot{X}(t)$ and, therefore,
$\bP(t)\dot{\cH}(t) = \mzero$. Applying Equation~(3.1) and Theorem~3.4, we conclude
that
$$
\mzero = \bP(t)\dot{\cH}(t) = -\bK(t)\cA(t) = -(1/2)\cA(t)\{\cA(t),t\}.
$$
The upshot is that the Schwarzian $\{\cA(t),t\}$ must be zero. \qed

As a corollary, we obtain an easy proof of the following important property of
the matrix Schwarzian (\cite[p. 205]{Zelikin:2000}).

\newclaim Corollary.
Let $M(t)$ be a smooth curve of $n \times n$ matrices such that $\dot{M}(t)$ is
invertible for all values of $t$. The matrix Schwarzian of $M(t)$ is identically
zero if and only if there exist matrices $A$, $B$, $C$, and $D$ such that
$M(t) = (C + tD)(A + tB)^{-1}$.

\Proof
Since the matrix Schwarzian of $M(t)$ is zero, so is the Schwarzian of the
fanning frame
$$
\cA(t) = \pmatrix{I \cr M(t) \cr} ,
$$
and of any other fanning frame spanning the same curve in the Grassmannian. In
particular, the Schwarzian of a normal frame $\cB(t) = \cA(t)X(t)$ vanishes
identically. By the preceding theorem, we have that
$\cB(t) = \cB(0) + t\dot{\cB}(0)$. If we write
$$
\cB(0) = \pmatrix{A \cr C \cr} \quad \hbox{\rm and}
\quad \dot{\cB}(0) = \pmatrix{B \cr D \cr},
\qquad \hbox{\rm then} \quad
\pmatrix{I \cr M(t) \cr} = \pmatrix{A + tB \cr C + tD \cr} X(t)^{-1}.
$$
It follows that $X(t) = (A + tB)$ and  $M(t) = (C + tD)(A + tB)^{-1}$. \qed

We now characterize the fanning curves of the form $\ell(t) = \exp(t\bX) \ell$.
These curves are ubiquitous in the study of the Jacobi equation in symmetric and
reductive spaces in Riemannian and Finsler geometry.

\newpar Definition.
A fanning curve $\ell(t)$ in the Grassmannian $\grassmann$ is said to be
{\sl parallel\/} if the Schwarzian of a normal frame spanning it is constant. It
is called {\sl weakly parallel\/} if whenever $\cA(t)$ is a normal frame
spanning $\ell(t)$, there exist $n \times n$ matrices $X$ and $Y$ such that
$\{\cA(t),t\} = \exp(-tY)\,X\,\exp(tY)$.

\newclaim Lemma.
A fanning frame $\cA(t)$ satisfies a constant-coefficient, second-order
differential equation $\ddot{\cA} + \dot{\cA}P + \cA Q = \mzero$ ($P$ and
$Q$ are $n \times n$ matrices) if and only if it is of the form
$\cA(t) = \exp(t\bX)\cA(0)$ for some linear transformation $\bX$ from $\R^{2n}$
to itself such that the $2n \times 2n$ matrix $\join{\cA(0)}{\bX\cA(0)}$ is
invertible.

\Proof
Assume that $\cA(t) = \exp(t\bX)\cA(0)$ for some $\bX$ satisfying the hypotheses
of the theorem. Since the columns of $\cA(0)$ and $\dot{\cA}(0) = \bX\cA(0)$
span $\R^{2n}$, there exist $n \times n$ matrices $P$ and $Q$ so that
$\bX^2\cA(0) + \bX\cA(0)P + \cA(0)Q = \mzero$. We have then that
$$
\ddot{\cA}(t) + \dot{\cA}(t)P + \cA(t)Q =
\exp(t\bX)\left(\bX^2\cA(0) + \bX\cA(0)P + \cA(0)Q \right) = \mzero .
$$
Conversely, we may write the equation
$\ddot{\cA}(t) + \dot{\cA}(t)P + \cA(t)Q = \mzero$ as
$$
\frac{d}{dt}\join{\cA(t)}{\dot{\cA}(t)} =
\join{\cA(t)}{\dot{\cA}(t)} \pmatrix{\mzero & -Q \cr I & -P} .
$$
It follows that $\cA(t) = \exp(t\bX)\cA(0)$, where
$$
\bX = \join{\cA(0)}{\dot{\cA}(0)} \pmatrix{\mzero & -Q \cr I & -P}
\join{\cA(0)}{\dot{\cA}(0)}^{-1} . \qed
$$

\newclaim Theorem.
A fanning curve in $\grassmann$ is weakly parallel if and only if it is of the
form $\ell(t) = \exp(t\bX)\ell(0)$, where $\bX$ is a linear transformation of
$\R^{2n}$ such that $\bX \ell(0)$ is transverse to $\ell(0)$.

\Proof
If the fanning curve $\ell(t)$ is weakly parallel and $\cA(t)$ is a normal frame
spanning it, there exist $n \times n$ matrices $X$ and $Y$ such that
$\ddot{\cA}(t) + (1/2)\cA(t)\exp(-tY)\,X\,\exp(tY) = \mzero$. Setting
$\cB(t) = \cA(t)\exp(-tY)$, we easily verify that $\cB(t)$ satisfies
the constant-coefficient, second-order differential equation
$$
\ddot{\cB}(t) + \dot{\cB}(t)(2Y) + \cB(t)(Y^2 + (1/2)X) = \mzero.
$$
By Lemma~5.4, this implies that $\cB(t) = \exp(t\bX)\cB(0)$ and, therefore,
$\ell(t) =  \exp(t\bX)\ell(0)$ for some linear transformation $\bX$ from
$\R^{2n}$ to itself.

Conversely, assume that $\ell(t) = \exp(t\bX)\ell(0)$ and let $\cA(0)$ be a
$2n \times n$ matrix whose columns span the subspace $\ell(0)$. By Lemma~5.4,
$\cA(t) = \exp(t\bX)\cA(0)$ satisfies a constant-coefficient, second-order
differential equation $\ddot{\cA} + \dot{\cA}P + \cA Q = \mzero$. Hence, the
Schwarzian $\{\cA(t),t\}$ is the constant matrix $X = 2Q - (1/2)P^2$, and the
Schwarzian of the normal frame $\cA(t)\exp(tP/2)$ is $\exp(-tP/2)\,X\,\exp(tP/2)$.
Thus, the fanning curve $\ell(t)$ is weakly parallel. \qed

\newclaim Corollary.
If $\cA(t)$ is any frame spanning a fanning curve $\ell(t)$ that is the projection
of a one-parameter subgroup of $\lgroup{2n}$ to the Grassmannian $\grassmann$,
then the Schwarzian of $\cA(t)$ satisfies a Lax equation
$$
\frac{d}{dt}\{\cA(t),t\} = [\{\cA(t),t\},Y(t)],
$$
where $Y(t)$ is a smooth curve of $n \times n$ matrices.

\Proof
Since $\ell(t)$ is weakly parallel, the Schwarzian of a normal frame $\cB(t)$
spanning it is of the form $\exp(-tY)X\exp(tY)$ for fixed $n \times n$
matrices $X$ and $Y$. Therefore, we have the equation
$$
\frac{d}{dt}\{\cB(t),t\} = [\{\cB(t),t\},Y].
$$
If $\cA(t)= \cB(t) X(t)$ is any other fanning frame spanning $\ell(t)$,
$$
\eqalign{
\frac{d}{dt}\, \{\cA(t),t\} &= \frac{d}{dt}\, X(t)^{-1} \{\cB(t),t\} X(t) \cr
                          &= -X(t)^{-1}\dot{X}(t) \{\cA(t),t\} +
                             X(t)^{-1}\frac{d}{dt} \{\cB(t),t\} X(t) +
                             \{\cA(t),t\}X(t)^{-1}\dot{X}(t) \cr
                          &= -X(t)^{-1}\dot{X}(t) \{\cA(t),t\} +
                             X(t)^{-1}[\{\cB(t),t\},Y] X(t) +
                             \{\cA(t),t\}X(t)^{-1}\dot{X}(t) \cr
                          &= [\{\cA(t),t\}, X(t)^{-1}YX(t) + X(t)^{-1}\dot{X}(t)].
                          \cr }
$$
\qed

\newclaim Theorem.
A fanning curve in $\grassmann$ is parallel if and only if it is of the
form $\ell(t) = \exp(t\bX)\ell(0)$, where $\bX$ is a linear transformation of
$\R^{2n}$ such that $\bX \ell(0)$ is transverse to $\ell(0)$ and
$\bX^2\ell(0) \subset \ell(0)$.

\Proof
If the fanning curve $\ell(t)$ is parallel, a normal frame $\cA(t)$ spanning it
satisfies a constant-coefficient, second-order differential equation
$\ddot{\cA} + \cA Q = \mzero$. By Lemma~5.4, this implies that
$\cA(t) = \exp(t\bX)\cA(0)$, where
$$
\bX = \join{\cA(0)}{\dot{\cA}(0)} \pmatrix{\mzero & -Q \cr I & \mzero}
\join{\cA(0)}{\dot{\cA}(0)}^{-1} .
$$
Notice that $\bX \cA(0) = \dot{\cA}(0)$ and $\bX^2 \cA(0) = -\cA(0)Q$. Therefore,
$\bX \ell(0)$ is transverse to $\ell(0)$ and $\bX^2 \ell(0) \subset \ell(0)$.

Conversely, let $\ell(t)$ be a fanning curve of the form $\exp(t\bX)\ell(0)$ for
some linear transformation $\bX$ from $\R^{2n}$ to itself. The condition that
$\bX^2\ell(0)$ be a subspace of  $\ell(0)$ implies that if $\cA(0)$ is any
$2n \times n$ matrix whose columns span $\ell(0)$, then $\cA(t) = \exp(t\bX)\cA(0)$
is a normal frame. By Lemma~5.4, the Schwarzian of $\cA(t)$ is constant and,
therefore, $\ell(t)$ is parallel. \qed

\newclaim Proposition.
The Jacobi endomorphism of a fanning curve $\ell(t)$ in $\grassmann$ is zero if
and only if it is of the form $\ell(t) = \exp(t\bX)\ell(0)$ with $\bX$ a linear
transformation $\bX$ from $\R^{2n}$ to itself such that $\bX\ell(0)$ is transversal
to $\ell(0)$ and $\bX^2$ is zero.

\Proof
If $\cA(0)$ is a $2n \times n$ matrix of whose columns span $\ell$, then
$\cA(t) = \exp(t\bX)\cA(0)$ is a normal frame spanning $\ell(t)$ and satisfying
$\ddot{\cA}(t) = \bX^2\cA(t) = \mzero$. It follows that the Jacobi endomorphism
of $\ell(t)$ is zero.

Conversely, if  $\cA(t)$ is a normal frame spanning a fanning curve with vanishing
Jacobi endomorphism, $\cA(t) = \cA(0) + t\dot{\cA}(0)$. If $\bX$ is any linear
transformation from $\R^{2n}$ to itself such that $\bX\cA(0) = \dot{\cA}(0)$ and
$\bX^2 = \mzero$, then $\cA(t) = \exp(t\bX)\cA(0)$ and, therefore,
$\ell(t) = \exp(t\bX)\ell(0)$.\qed

\mysection{Fanning curves of Lagrangian subspaces}
We now turn to the study of fanning curves of Lagrangian subspaces on $\R^{2n}$
provided with its standard symplectic structure
$$
\omega(\bv,\bw) = \bv^\trans \bJ \bw, \quad \hbox{\rm where} \quad
\bJ = \pmatrix{\mzero & -I \cr I & \mzero}.
$$

Recall that the set $\Lambda_n$ of all Lagrangian subspaces of $\R^{2n}$---the
{\sl Lagrangian Grassmannian}---is a smooth submanifold of the Grassmannian
$\grassmann$. If $\ell$ is a Lagrangian subspace of $\R^{2n}$ the tangent space
$T_{\ell}\Lambda_n$ is canonically isomorphic to the space of symmetric bilinear
forms on $\ell$. It will be useful to describe this isomorphism in terms of frames.

Let us say that a {\sl Lagrangian frame\/} is a $2n \times n$ matrix $\cA$ of rank
$n$ such that the subspace spanned by the columns of $\cA$ is Lagrangian.
Equivalently, a $2n \times n$ matrix $\cA$ is a Lagrangian frame if it has rank $n$
and $\cA^\trans \bJ \cA = \mzero$. If $\cA(t)$ is a smooth curve of Lagrangian
frames spanning a curve $\ell(t)$ of Lagrangian subspaces, we have that
$$
\cA(t)^\trans \bJ \cA(t) = \mzero \quad \hbox{\rm and} \quad
\dot{\cA}(t)^\trans \bJ \cA(t) + \cA(t)^\trans \bJ \dot{\cA}(t) = \mzero.
\leqno{(6.1)}
$$
The second equation is simply the derivative of the first.

\newpar Definition.
The {\sl Wronskian\/} of a smooth curve $\cA(t)$ of Lagrangian frames is the
curve of $n \times n$ matrices defined by $W(t) = -\cA(t)^\trans \bJ \dot{\cA}(t)$.

Note that in two dimensions the Wronskian of a curve of Lagrangian frames
$(f_1(t), f_2(t))^{\trans}$ is the ``usual'' Wronskian
$f_1(t)\dot{f}_2(t) - f_2(t)\dot{f}_1(t)$.

\newclaim Proposition.
The Wronskian $W(t)$ of a smooth curve $\cA(t)$ of Lagrangian frames satisfies the
following properties:
\item{1.} $W(t)$ is symmetric for all values of $t$.
\item{2.} If $X(t)$ is a smooth curve of invertible $n \times n$ matrices, the
           Wronskian of the curve of Lagrangian frames $\cA(t)X(t)$ is
           $X(t)^\trans W(t) X(t)$.
\item{3.} If $\bS$ is a linear symplectic transformation of ($\R^{2n},\omega)$,
           the Wronskian of $\bS \cA(t)$ is also $W(t)$.
\item{4.} The curve of Lagrangian frames $\cA(t)$ is fanning if and only if
           $W(t)$ is invertible for all values of $t$.
\item{5.} If $s$ is a diffeomorphism of the real line, the Wronskian of $\cA(s(t))$
           is $W(s(t))\dot{s}(t)$.

\Proof
Properties~(1) and~(2) follow immediately from Equation~(6.1) and the definition
of the Wronskian. Property~(3) follows from the fact that $\bS$ is symplectic
(i.e., $\bS^\trans \bJ \bS = \bJ$), while Property~(5) is a trivial computation.
Lastly, to show that $\cA(t)$ is fanning if and only if $W(t)$ is never singular,
it is enough to remark that $\join{\cA(t)}{\dot{\cA}(t)}$ is nonsingular if and
only if
$$
\join{\cA(t)}{\dot{\cA}(t)}^\trans \, \bJ \, \join{\cA(t)}{\dot{\cA}(t)} =
\pmatrix{ \mzero & -W(t) \cr W(t) & \dot{\cA}(t)^\trans \bJ \dot{\cA}(t) \cr}
$$
is nonsingular.
\qed

The geometric interpretation of Property~(2) is that if $\ell(t)$ is a smooth
curve of Lagrangian subspace of $\R^{2n}$ and $\cA(t)$ is a curve of Lagrangian
frames spanning it, the bilinear form on the subspace $\ell(\tau)$ defined by the
matrix $W(\tau)$ in the basis formed by the columns of $\cA(\tau)$ is well defined
and does not depend on the curve of frames spanning $\ell(t)$. This gives the
canonical isomorphism between $T_{\ell(\tau)}\Lambda_n$ and the space of bilinear
forms on $\ell(\tau)$.

Properties~(1) and ~(4) imply that the Wronskian $W(t)$ of a fanning curve of
Lagrangian frames is a curve of invertible, symmetric matrices and, therefore,
that the index of the matrices $W(t)$ is constant. By the remark in the preceding
paragraph, the index of $W(t)$ depends only on the fanning curve on the Lagrangian
Grassmannian and not on the particular fanning frame that spans it.

\newpar Definition.
The {\sl signature\/} of a fanning curve of Lagrangian subspaces is the
index of the  Wronskian of any curve of Lagrangian frames that spans it.

Before considering the congruence problem for fanning curves of Lagrangian
subspaces, let us introduce a large class of examples that arise in applications
and whose invariants are easily computed.

\subsection{Systems of Lagrange equations}
Let us consider a system of Lagrange equations of form
$$
\frac{d}{dt}(\dot{\bx}^\trans K(t)) + \bx^\trans V(t) = \bzero, \leqno{(6.2)}
$$
where $K(t)$ and $V(t)$ are smooth curves of $n \times n$ symmetric matrices
and $K(t)$ is invertible for all values of $t$. The vectors
$\bx^\trans$ and $\dot{\bx}^\trans$ in $\R^n$ are written as $1 \times n$ matrices.
We may also write Equation~(6.2) in Hamiltonian form
$$
\dot{\bx}^\trans = \bp^\trans K(t)^{-1} \, , \qquad
\dot{\bp}^\trans = - \bx^\trans V(t)
$$
with Hamiltonian
$H(\bx,\bp,t) = (1/2) \left(\bp^\trans K(t)^{-1} \bp + \bx^\trans V(t) \bx \right)$.

\newclaim Proposition.
If $\cA(t)$ is a $2n \times n$ matrix whose rows are solutions of Equation~(6.2)
and such that $\join{\cA(0)}{\dot{\cA}(0)K(0)}$ is a symplectic matrix, then
$\cA(t)$ is a fanning curve of Lagrangian frames whose Wronskian is
$W(t) = K(t)^{-1}$ and whose Schwarzian is
$$
\{\cA(t),t\} = 2V(t)K(t)^{-1} - (1/2)\left(\dot{K}(t)K(t)^{-1}\right)^2 -
\frac{d}{dt} \left(\dot{K}(t)K(t)^{-1}\right).
$$

\Proof
Notice that Equation~(6.2) may be written as
$$
\frac{d}{dt} \join{\cA(t)}{\dot{\cA}(t)K(t)} =
\join{\cA(t)}{\dot{\cA}(t)K(t)}
\pmatrix{\mzero & -V(t) \cr K(t)^{-1} & \mzero}
$$
and that, since both $K(t)$ and $V(t)$ are symmetric, the last matrix defines a
curve in the Lie algebra of the group of symplectic transformations of
$(\R^{2n}, \omega)$. It follows that if the initial condition
$\join{\cA(0)}{\dot{\cA}(0)K(0)}$ is symplectic, then
$\join{\cA(t)}{\dot{\cA}(t)K(t)}$ is symplectic for all values of $t$. This
implies at once that $\cA(t)$ is a curve of Lagrangian frames. Moreover, since
$\join{\cA(t)}{\dot{\cA}(t)}$ is nonsingular if and only if
$\join{\cA(t)}{\dot{\cA}(t)K(t)}$ is nonsingular, we have that $\cA(t)$ is a
fanning curve of frames.

Note that since $\join{\cA(t)}{\dot{\cA}(t)K(t)}$ is symplectic,
$\cA(t)^\trans \bJ \dot{\cA}(t)K(t) = -I$ and, therefore, the Wronskian
$W(t) = -\cA(t)^\trans \bJ \dot{\cA}(t) = K(t)^{-1}$. The Schwarzian of $\cA(t)$
is easily computed using Equations~(1.2) and~(6.2).
\qed

\subsection{Lagrangian normal frames and solution of the congruence problem}
Let us start with a simple---albeit useful---remark:

\newclaim Lemma.
If $\cA(t)$ is a curve of Lagrangian frames such the columns of $\ddot{\cA}(t)$
are linear combinations of those of $\cA(t)$, then
$\dot{\cA}(t)^\trans \bJ \dot{\cA}(t) = \mzero$. In particular, the Wronskian
of such a frame is constant.

\Proof
Differentiating the equation
$\dot{\cA}(t)^\trans \bJ \cA(t) + \cA(t)^\trans \bJ \dot{\cA}(t) = \mzero$ and
using that $\ddot{\cA}(t)^\trans \bJ \cA(t)$ and $\cA(t)^\trans \bJ \ddot{\cA}(t)$
are both zero, we conclude that
$$
\mzero = \ddot{\cA}(t)^\trans \bJ \cA(t) + 2 \dot{\cA}(t)^\trans \bJ \dot{\cA}(t)
          + \cA(t)^\trans \bJ \ddot{\cA}(t)
       = 2 \dot{\cA}(t)^\trans \bJ \dot{\cA}(t).
$$

Since the derivative of the Wronskian of $\cA(t)$ is
$-\dot{\cA}(t)^\trans \bJ \dot{\cA}(t) - \cA(t)^\trans \bJ \ddot{\cA}(t)$
and both terms are zero, the Wronskian is constant.
\qed

In view of this lemma, if a Lagrangian frame $\cA(t)$ is normal in the sense
of Section~4, then its Wronskian is constant. This allows us to define a more
restrictive class of normal frames in the Lagrangian setting.

\newpar Definition.
A fanning Lagrangian frame $\cA(t)$ is said to be a {\sl Lagrangian
normal frame\/} if the columns of $\ddot{\cA}(t)$ are linear combinations of
those of $\cA(t)$ and its Wronskian is the diagonal $n \times n$ matrix
$I_{n,k}$ whose first $k$ $(0 \leq k \leq n)$ diagonal entries equal $-1$ and the
remaining diagonal entries equal $1$.

\newclaim Proposition.
If $\ell(t)$ is a fanning curve of Lagrangian subspaces in $\R^{2n}$, there exists
a Lagrangian normal frame that spans it. Moreover if the signature of $\ell(t)$ is
$k$ $(0 \leq k \leq n)$ and $\cA(t)$ and $\cB(t)$ are two Lagrangian normal frames
spanning $\ell(t)$, there exists a fixed invertible $n \times n$ matrix $X$ that
preserves the quadratic form defined by the matrix $I_{n,k}$ and such that
$\cB(t) = \cA(t)X$.

\Proof
By Proposition~4.2, there exists a normal frame $\cA(t)$ spanning $\ell(t)$ and,
by Lemma~6.5, the Wronskian of $\cA(t)$ is a fixed symmetric matrix $W$. If
$X$ is any $n \times n$ matrix such that $X^\trans W X = I_{n,k}$, where
$k$ is the index of $W$, then $\cA(t)X$ is a Lagrangian normal frame spanning
$\ell(t)$.

Two Lagrangian normal frames $\cA(t)$ and $\cB(t)$ spanning the same curve $\ell(t)$
on the Lagrangian Grassmannian are, in particular, two normal frames spanning the
same curve of $n$-dimensional subspaces of $\R^{2n}$. It follows that there exists
an invertible $n \times n$ matrix $X$ such that $\cB(t) = \cA(t)X$. Since the
Wronskians of $\cA(t)$ and $\cB(t)$ are both equal to $I_{n,k}$, Property~(2) in
Proposition~6.2 implies that $X^\trans I_{n,k} X = I_{n,k}$.
\qed

Lagrangian normal frames spanning a given fanning curve of Lagrangian
subspaces are therefore determined up to multiplication by a fixed
matrix in the orthogonal group $O(n-k,k)$.

\newclaim Theorem.
Two fanning curves of Lagrangian subspaces of $\R^{2n}$ are congruent if
and only if they have the same signature $k$ $(0 \leq k \leq n)$ and the
Schwarzians of any two of their Lagrangian normal frames are conjugate by a
constant $n \times n$ matrix in $O(n-k,k)$.

\Proof
The proof of this theorem is similar to that of Theorem~4.3, except in one
important point: having reduced the proof to verifying that two Lagrangian normal
frames $\cA(t)$ and $\cB(t)$ that have the same Wronskian and Schwarzian are
congruent, it is somewhat more tricky to find a linear symplectic transformation
such that $\bS\cA(0) = \cB(0)$ and $\bS \dot{\cA}(0) = \dot{\cB}(0)$.
Just as in the proof of Theorem~4.3, this would imply that $\bS \cA(t)$ equals
$\cB(t)$---and prove the theorem---because both Lagrangian normal frames satisfy the
same second-order differential equation and have the same initial conditions.

In order to find $\bS$ we remark that the matrices
$\bA := \join{\cA(0)}{\dot{\cA}(0)I_{n,k}}$ and
$\bB := \join{\cB(0)}{\dot{\cB}(0)I_{n,k}}$ are symplectic. Indeed, using that
$\cA(0)^\trans \bJ \cA(0) = \mzero$,
$\dot{\cA}(0)^\trans \bJ \dot{\cA}(0) = \mzero$ (Lemma~6.5), and that both
$- \cA(0)^\trans \bJ \dot{\cA}(0)$ and  $\dot{\cA}(0)^\trans \bJ \cA(0)$ are
equal to the Wronskian $I_{n,k}$, we have that
$$
\join{\cA(0)}{\dot{\cA}(0)I_{n,k}}^\trans \, \bJ \,
\join{\cA(0)}{\dot{\cA}(0)I_{n,k}}
= \pmatrix{ \mzero & -I_{n,k}^2 \cr I_{n,k}^2 & \mzero} = \bJ .
$$

The matrix $\bS := \bB \bA^{-1}$ is the one we search.
\qed

The use of normal Lagrangian frames also makes it easy to prove that the
horizontal curve of a fanning curve of Lagrangian subspaces is itself a
curve of Lagrangian subspaces, and that the Schwarzian of a fanning
Lagrangian frame $\cA(t)$ is symmetric with respect to the inverse of
its Wronskian (i.e. that the matrix $\{\cA(t),t\}W(t)^{-1}$ is symmetric).

\newclaim Proposition.
The horizontal curve of a fanning curve of Lagrangian subspaces is itself a
curve of Lagrangian subspaces. Equivalently, if $\cA(t)$ is a fanning Lagrangian
frame, then its horizontal derivative $\cH(t)$ is a curve of Lagrangian frames.

\Proof
Recall that if  $\cB(t)$ is normal frame, then its horizontal derivative is
simply $\dot{\cB}(t)$. If $\cB(t)$ is also a Lagrangian normal frame, then
Lemma~6.5 implies that $\dot{\cB}(t)$ is a curve of Lagrangian frames.

If $\cA(t)$ is any curve of Lagrangian frames, there exist a normal Lagrangian
frame $\cB(t)$ and  a curve of invertible $n \times n$ matrices so that
$\cA(t)= \cB(t)X(t)$. By Proposition~2.9, the horizontal derivative of
$\cA(t)$ is $\dot{\cB}(t)X(t)$, which is a curve of Lagrangian frames.
\qed

\newclaim Proposition.
If $W(t)$ denotes the Wronskian of a fanning Lagrangian frame $\cA(t)$,
the matrix $\{\cA(t),t\}W(t)^{-1}$ is symmetric.

\Proof
Let us first assume that $\cB(t)$ is a normal Lagrangian frame with Wronskian
$I_{n,k}$. An argument identical to the last step in the proof of Theorem~6.8
shows that the matrix $\join{\cB(t)}{\dot{\cB(t)}I_{n,k}}$ is symplectic. Since
$\cB(t)$ is a normal frame, we have that $\ddot{\cB} + (1/2)\cB(t)\{\cB(t),t\}$
and, therefore,
$$
\frac{d}{dt} \join{\cB(t)}{\dot{\cB(t)}I_{n,k}} =
\join{\cB(t)}{\dot{\cB(t)}I_{n,k}}
\pmatrix{\mzero & -\{\cB(t),t\}I_{n,k} \cr I_{n,k} & \mzero} .
$$
It follows that the rightmost  matrix must be in the Lie algebra of the group
of linear symplectic transformations and, consequently, that $\{\cB(t),t\}I_{n,k}$
is symmetric.

If $\cA(t)$ is any fanning Lagrangian frame, there exits a Lagrangian normal frame
$\cB(t)$ and a curve of $n \times n$ invertible matrices $X(t)$ such that
$\cA(t) = \cB(t)X(t)$. It follows that the Wronskian of $\cA(t)$ is
$W(t) = X(t)^\trans I_{n,k}X(t)$ and that $\{\cA(t),t\} = X^{-1}\{\cB(t),t\}X(t)$.
Therefore, the matrix
$\{\cA(t),t\} W(t)^{-1} = X^{-1}\{\cB(t),t\}I_{n,k}(X^{-1})^\trans$ is symmetric.
\qed

\mysection{Geometry of one- and two-jets of fanning curves}
In order to understand the geometry of submanifolds in a homogeneous space, it is
necessary to study the prolonged action of the group of symmetries on the
spaces of jets. Except in simple classical cases, this is usually an arduous
undertaking; specially in the case of noncompact symmetry groups. It is
then a pleasant surprise that the action of $\lgroup{2n}$ on the spaces of jets
of fanning curves in $\grassmann$ can be easily described and that much geometric
information and insight can be gained from its study.

\subsection{One-jets of fanning curves and the fundamental endomorphism}
The easiest way to describe the action of $\lgroup{2n}$ on the space
$\kfan{k}{\grassmann}$ of  $k$-jets of fanning curves on the Grassmannian
is to use the natural representation of $\kfan{k}{\grassmann}$ as the quotient
of the space $\kfan{k}{\frames}$ of $k$-jets of fanning frames by the action of
the group $\kmaps{k}$ of $k$-jets of smooth curves of invertible $n \times n$
matrices. For instance, the space of one-jets of fanning frames consists of ordered
pairs $(\cA, \dot{\cA})$ of $2n \times n$ matrices such that
$\join{\cA}{\dot{\cA}}$ is invertible, and the group of one-jets of curves of
invertible $n \times n$ matrices consists of ordered pairs $(X,\dot{X})$ of
$n \times n$ matrices where $X$ is invertible. The actions of $\kmaps{1}$ and
$\lgroup{2n}$ on $\kfan{1}{\frames}$ are given by
$$
(\cA,\dot{\cA}) \cdot (X, \dot{X}) = (\cA X, \dot{\cA} X + \cA \dot{X})
\quad \hbox{\rm and} \quad
\bT \cdot (\cA, \dot{\cA}) = (\bT\cA, \bT\dot{\cA}) ,
$$
respectively. If we represent the elements of
$\kfan{1}{\grassmann} = \kfan{1}{\frames}/\kmaps{1}$ as equivalence classes
$[(\cA,\dot{\cA})]$, the action of $\lgroup{2n}$ on the one-jets of fanning
curves in the Grassmannian is simply described by
$\bT \cdot [(\cA, \dot{\cA})] = [(\bT\cA, \bT\dot{\cA})]$.

\newclaim Proposition.
The group of invertible linear transformations of $\R^{2n}$ acts transitively
on the space of one-jets of fanning frames and, a fortiori, on the space of
one-jets of fanning curves in the Grassmannian $\grassmann$.

\Proof
If $(\cA, \dot{\cA}) \in \kfan{1}{\frames}$, then $\bT := \join{\cA}{\dot{\cA}}$
is in $\lgroup{2n}$ and
$$
\bT^{-1} \cdot (\cA, \dot{\cA}) =
\left( \pmatrix{I \cr \mzero}, \pmatrix{\mzero \cr I} \right) . \qed
$$

\newclaim Theorem.
A map from the space of one-jets of fanning curves in $\grassmann$ to the Lie
algebra $\lalgebra$ is equivariant with respect to the $\lgroup{2n}$ action on
these spaces if and only if it is of the form $a\bI + b\bF$, where $\bI$ is the
identity matrix and $a$ and $b$ are real numbers.

\Proof
Let $\bG : \kfan{1}{\frames} \rightarrow \lalgebra$ be a map that is
invariant under the action of $\kmaps{1}$ and equivariant with respect to
the action of $\lgroup{2n}$. Writing $\bG(\cA,\dot{\cA})$ as
$$
\bG(\cA,\dot{\cA}) =
\join{\cA}{\dot{\cA}}
\pmatrix{G_{11}(\cA,\dot{\cA}) & G_{12}(\cA,\dot{\cA}) \cr
         G_{21}(\cA,\dot{\cA}) & G_{22}(\cA,\dot{\cA})}
\join{\cA}{\dot{\cA}}^{-1}
$$
and using that $\bG(\bT\cA,\bT\dot{\cA}) = \bT \bG(\cA,\dot{\cA}) \bT^{-1}$,
we obtain that
$$
\pmatrix{G_{11}(\bT\cA,\bT\dot{\cA}) & G_{12}(\bT\cA,\bT\dot{\cA}) \cr
         G_{21}(\bT\cA,\bT\dot{\cA}) & G_{22}(\bT\cA,\bT\dot{\cA})} =
\pmatrix{G_{11}(\cA,\dot{\cA}) & G_{12}(\cA,\dot{\cA}) \cr
         G_{21}(\cA,\dot{\cA}) & G_{22}(\cA,\dot{\cA})}
$$
for all matrices $\bT \in \lgroup{2n}$. Since $\lgroup{2n}$ acts transitively on
the space of one-jets of fanning frames, this implies that the blocks $G_{ij}$ are
constant $n \times n$ matrices. Moreover, the invariance under the action of
$\kmaps{1}$ imposes the condition
$$
\pmatrix{X & \dot{X} \cr \mzero & X}
\pmatrix{G_{11} & G_{12} \cr
         G_{21} & G_{22}}
\pmatrix{X & \dot{X} \cr \mzero & X}^{-1}
= \pmatrix{G_{11} & G_{12} \cr
         G_{21} & G_{22}}
$$
for all choices of matrices $X$ and $\dot{X}$ with $X$ invertible. This easily
implies that
$$
\pmatrix{G_{11} & G_{12} \cr G_{21} & G_{22}}
= \pmatrix{aI & bI \cr \mzero & aI}
$$
for some real numbers $a$ and $b$ and, therefore, that $\bG = a\bI + b \bF$. \qed

\subsection{Two-jets of fanning curves and the horizontal map}
The space of two-jets of fanning frames is the set of ordered triples
$(\cA,\dot{\cA},\ddot{\cA})$ of $2n \times n$ matrices such that
$\join{\cA}{\dot{\cA}}$ is invertible. Likewise, the two-jets of
curves of $n \times n$ invertible matrices are ordered triples of $n \times n$
matrices $(X, \dot{X}, \ddot{X})$ such that $X$ is invertible. The actions of
$\kmaps{2}$ and $\lgroup{2n}$ on $\kfan{2}{\frames}$ are given by
$$
(\cA,\dot{\cA},\ddot{\cA}) \cdot (X, \dot{X}, \ddot{X}) =
(\cA X, \dot{\cA} X + \cA \dot{X}, \ddot{\cA}X + 2\dot{\cA}\dot{X} + \cA\ddot{X})
\quad \hbox{\rm and} \quad
\bT \cdot (\cA, \dot{\cA},\ddot{\cA}) = (\bT\cA, \bT\dot{\cA},\bT\ddot{\cA}) ,
$$
respectively. However, unlike the case of one-jets, the action of $\lgroup{2n}$
on the space of two-jets of fanning frames is not transitive. Indeed, it is easy
to show that the two-jets $(\cA,\dot{\cA},\ddot{\cA})$ and
$(\cB,\dot{\cB},\ddot{\cB})$ are in the same $\lgroup{2n}$-orbit if and only if
the matrices $\join{\cA}{\dot{\cA}}^{-1}\ddot{\cA}$ and
$\join{\cB}{\dot{\cB}}^{-1} \ddot{\cB}$ are equal. Nevertheless, we still have that
$\lgroup{2n}$  acts transitively on $\kfan{2}{\grassmann}$.

\newclaim Proposition.
The group of invertible linear transformations of $\R^{2n}$ acts transitively
on the space of two-jets of fanning curves in $\grassmann$.

\Proof
Since $\kfan{2}{\grassmann} = \kfan{2}{\frames}/\kmaps{2}$, all that needs to
be shown is that the joint action of $\lgroup{2n}$ and $\kmaps{2}$ on the
space of two-jets of fanning frames is transitive. To verify this, just note that
if $(\cA,\dot{\cA},\ddot{\cA}) \in \kfan{2}{\frames}$ and
$\ddot{\cA} + \dot{\cA}P + \cA Q = \mzero$, then
$$
\join{\cA}{\dot{\cA} + (1/2)\cA P}^{-1} \cdot (\cA,\dot{\cA},\ddot{\cA})
\cdot (I, (1/2)P,Q) =
\left(
\pmatrix{I \cr \mzero}, \pmatrix{\mzero \cr I}, \pmatrix{\mzero \cr \mzero}
\right). \qed
$$

We are now ready to prove the characterization of the horizontal curve of
a fanning curve in the Grassmannian $\grassmann$ that we stated in the
introduction.

\newclaim Theorem.
The assignment that sends a fanning curve $\ell(t)$ to its horizontal curve
$h(t)$ is characterized by the following four properties:
\item{1.} At each time $t$ the subspace $h(t)$ is transversal to $\ell(t)$;
\item{2.} The subspace $h(\tau)$ depends only on the two-jet of the curve $\ell(t)$
          at $t = \tau$.
\item{3.} If $\bT$ is an invertible linear transformation of $\R^{2n}$, the
          horizontal curve of $\bT\ell(t)$ is $\bT h(t)$.
\item{4.} If $\ell(t)$ is spanned by a line $\cA + t\cB$ in the space of frames,
          $h(t)$ is constant.

\newclaim Lemma.
A  $\lgroup{2n}$-equivariant map $j : \kfan{2}{\grassmann} \rightarrow \grassmann$
is such that the subspaces $j([(\cA,\dot{\cA},\ddot{\cA})])$ and $[\cA]$ are
always transversal if and only if it is of the form
$$
[(\cA,\dot{\cA},\ddot{\cA})] \longmapsto [\dot{\cA} + (1/2)\cA P + a\cA],
\leqno{(7.1)}
$$
where $a$ is some real number and $P$ is the $n \times n$ matrix defined by the
equation $\ddot{\cA} + \dot{\cA}P + \cA Q = \mzero$.

\Proof
Notice that, at the level of two-jets of fanning frames,
$(\cA,\dot{\cA},\ddot{\cA}) \mapsto \dot{\cA} + (1/2)\cA P$ is just
the horizontal derivative we studied in Section~2. For any real number $a$,
the subspace $[\dot{\cA} + (1/2)\cA P + a\cA]$ is transversal to $[\cA]$ and,
using the properties of the horizontal derivative, it is easy to show that the
map defined by Equation~(7.1) is $\lgroup{2n}$ equivariant.

In order to prove the converse, let us define
$\bP : \kfan{2}{\frames} \rightarrow \lalgebra$ as the map whose value
at a two-jet $(\cA,\dot{\cA},\ddot{\cA}) \in \kfan{2}{\frames}$ is the
projection with range $[\cA]$ and kernel $j([(\cA,\dot{\cA},\ddot{\cA})])$.
The properties of $j$ are equivalent to the following properties of $\bP$:
$$
\eqalign{
& \hbox{\rm (i)}\,
\bP(\cA,\dot{\cA},\ddot{\cA})^2 = \bP(\cA,\dot{\cA},\ddot{\cA}), \cr
& \hbox{\rm (iii)}\,
\bP(\bT\cA,\bT\dot{\cA},\bT\ddot{\cA}) =
                       \bT \bP(\cA,\dot{\cA},\ddot{\cA}) \bT^{-1}}
\qquad
\eqalign{
&\hbox{\rm (ii)}\,
\bP(\cA,\dot{\cA},\ddot{\cA})\cA = \cA, \cr
&\hbox{\rm (iv)}\,
\bP(\cA X, \dot{\cA} X + \cA \dot{X}, \ddot{\cA}X + 2\dot{\cA}\dot{X}
              + \cA\ddot{X}) = \bP(\cA,\dot{\cA},\ddot{\cA}).}
$$

Using (i) and (ii), we see that there exists a function
$R = R(\cA,\dot{\cA},\ddot{\cA})$ with values in the space of $n \times n$
matrices such that
$$
\bP(\cA,\dot{\cA},\ddot{\cA}) =
\join{\cA}{\dot{\cA}}
\pmatrix{I & R \cr
         \mzero & \mzero}
\join{\cA}{\dot{\cA}}^{-1} .
$$
Moreover, since $\join{\cA}{\dot{\cA}}^{-1}\ddot{\cA}$ is the complete invariant
for the action of $\lgroup{2n}$ on two-jets of fanning frames, Property (iii)
implies that $R(\cA,\dot{\cA},\ddot{\cA})$ depends only on
$\join{\cA}{\dot{\cA}}^{-1}\ddot{\cA}$. With a slight abuse of notation, we
write $R = R(\join{\cA}{\dot{\cA}}^{-1}\ddot{\cA})$.

We now claim that (iv) is equivalent to the condition that
$R(\join{\cA}{\dot{\cA}}^{-1}\ddot{\cA})$ is of the form $ -(1/2)P - aI$, where
$a$ is an arbitrary real number. This is a somewhat tedious verification whose only
underlying idea is to break up Property~(iv) into three simpler invariance
properties: whenever $X$, $\dot{X}$,
and $\ddot{X}$ are $n \times n$ matrices with $X$ invertible,
$$
\eqalign{
&\hbox{\rm (a)}\,
\bP(\cA,\dot{\cA},\ddot{\cA} + \cA\ddot{X}) = \bP(\cA,\dot{\cA},\ddot{\cA}),
\quad
\hbox{\rm (b)}\,
\bP(\cA X,\dot{\cA}X,\ddot{\cA}X) = \bP(\cA,\dot{\cA},\ddot{\cA}), \cr
& \hbox{\rm (c)}\,
\bP(\cA,\dot{\cA} + \cA\dot{X},\ddot{\cA} + 2\dot{\cA}\dot{X}) =
              \bP(\cA,\dot{\cA},\ddot{\cA}).}
$$

From (a) we obtain that $R(\join{\cA}{\dot{\cA}}^{-1}\ddot{\cA})$ is really just a
function of the last $n$ rows of its argument. In other words, with a some abuse
of notation, we may write $R(\join{\cA}{\dot{\cA}}^{-1}\ddot{\cA}) = R(-P)$. From
(b) we obtain that $R(-X^{-1}PX) = X^{-1}R(-P)X$ for all $n \times n$ invertible
matrices $X$. In particular, this implies that $R(\mzero) = -aI$ for some real
number $a$. From (c) we obtain that $R(-P + 2 \dot{X}) = R(-P) + \dot{X}$. On
setting $2\dot{X} = P$, we have that
$R(-P) = -(1/2)P + R(\mzero) = -(1/2)P - aI$.

Since $j([(\cA,\dot{\cA},\ddot{\cA})])$ is the kernel of
$$
\bP(\cA,\dot{\cA},\ddot{\cA}) =
\join{\cA}{\dot{\cA}}
\pmatrix{I &   -(1/2)P - aI \cr
         \mzero & \mzero}
\join{\cA}{\dot{\cA}}^{-1},
$$
it must equal $[\dot{\cA} + (1/2)\cA P + a\cA]$.
\qed

\noindent
{\sl Proof of Theorem~7.4.\/}
Using the Properties~1--3 of $h(t)$, it follows from the previous lemma that
$h(t) = [\dot{\cA}(t) + (1/2)\cA(t)P(t) + a \cA(t)]$ for any choice of frame
$\cA(t)$ spanning $\ell(t)$.  When $\cA(t)$ is a line $\cA + t\cB$ in the space
of frames, then $h(t) = [\cB + a(\cA + t\cB)]$ is constant if and only if $a$ is
zero.
\qed
\mysection{A short survey of other approaches}
In this final section we briefly consider other approaches to the geometry of
fanning curves and exhibit their relationship to the invariants studied in
this paper. The reader is referred to the original papers for many of the proofs
and technical details.

\subsection{S. Ahdout's definition of the horizontal curve}
In the interesting paper~\cite{Ahdout:1989}, S.~Ahdout studies fanning curves
of Lagrangian subspaces and submanifolds and, as an application, shows that a
convex body in Euclidean $n$-space is determined up to rigid motion by the
billiard system it defines.  Most of the paper is concerned with the study of
fanning curves of Lagrangian submanifolds through the use of normal forms, but
it is there that fanning curves of Lagrangian subspaces and their horizontal
curves explicitly appear in the literature for the first time.

To a triple of $n$-dimensional subspaces $(\ell_0,\ell;\ell_\infty)$ in $\R^{2n}$
such that $\ell_\infty$ is transverse to both $\ell_0$ and $\ell$, we assign the
linear transformation
$\phi(\ell_0,\ell;\ell_\infty) : \ell_0 \rightarrow \R^{2n}/\ell_0$
that equals the composition of canonical projection from $\R^{2n}$ to the
quotient space $\R^{2n}/\ell_0$ and the unique linear transformation
from $\ell_0$ to $\ell_\infty$ whose graph in $\ell_0 \oplus \ell_\infty = \R^{2n}$
equals $\ell$. In what follows $\phi$ is considered as a function of the variable
$\ell$ whose domain is the open subset of $\grassmann$ consisting of all
$n$-dimensional subspaces transversal to $\ell_0$. This last subspace and
$\ell_\infty$ are treated as a parameters. It is well-known that if $\ell_\infty'$
is transversal to both $\ell_0$ and $\ell$, the differentials of
$\phi(\ell_0,\cdot\,;\ell_\infty)$ and $\phi(\ell_0,\cdot\,;\ell_\infty')$ at
$\ell$ are equal and describe the canonical isomorphism between
$T_{\ell_0}\grassmann$ and $\hom{\ell_0}{\R^{2n}/\ell_0}$.

A curve $\ell(t)$ in $\grassmann$ is fanning if at every time $\tau$ the linear map
$\frac{d}{dt}\phi(\ell(\tau),\ell(t);\ell_\infty)$ is invertible at $t = \tau$.
Here $\ell_\infty$ is any $n$-dimensional subspace transversal to $\ell(\tau)$
and, therefore, transversal to $\ell(t)$ for $t$ close to $\tau$.
In Theorem~(1.3) of \cite{Ahdout:1989}, Ahdout proves that {\sl if $\ell(t)$ is a
fanning curve in $\grassmann$, at every fixed time $\tau$ there exists a unique
$n$-dimensional subspace $h(\tau)$ that is transversal to $\ell(\tau)$ and such
that
$$
\frac{d^2}{dt^2}\phi(\ell(\tau),\ell(t);h(\tau))
$$
is zero at $t = \tau$.} A simple application of Theorem~7.4 shows that this
is indeed another definition of the horizontal curve.

\subsection{The approach of A.~A.~Agrachev et al. to the geometry of fanning curves}
In the series of papers \cite{Agrachev-Gamkrelidze:1997}, \cite{Agrachev:1998},
\cite{Agrachev-Zelenko:2001}, \cite{Agrachev-Zelenko:2002a},
and \cite{Agrachev-Zelenko:2002b}, A.~A.~Agrachev, R.~Gamkrelidze, and I.~Zelenko
study the geometry of curves in the Lagrangian Grassmannian by means of an
ingenious approach based on the use of Laurent series of matrices. The following
interpretation of their approach to the study of fanning curves touches only a
fraction of their work, most of which is concerned with special classes of
non-fanning curves arising from applications in control theory. We mention
in passing that Agrachev  uses the term ``regular Jacobi curves" to refer to
fanning curves in the Lagrangian Grassmannian and the term ``derivative curve''
to refer to the horizontal curve.

Let $(\ell_0;\ell_\infty,\ell)$ be a triple of $n$-dimensional subspaces in
$\R^{2n}$  such that $\ell_\infty$ and $\ell$ are transverse to $\ell_0$, and let
$A$ be the linear transformation from $\ell_\infty$ to $\ell_0$ whose graph in
$\ell_\infty \oplus \ell_0 = \R^{2n}$ is $\ell$. The nilpotent linear
transformation $\bN(\ell_0;\ell_\infty,\ell)$ from $\R^{2n}$ to itself defined
by the fact that its restriction to $\ell_0$ is identically zero and its
restriction to $\ell_\infty$ coincides with $A$ is characterized by three properties:
(1) its kernel contains $\ell_0$, (2) its range is contained in $\ell_0$, and
(3) the image of $\ell_\infty$ under the transformation
$\bI + \bN(\ell_0;\ell_\infty,\ell)$ is equal to $\ell$.
Central to Agrachev's approach are the two following properties of $\bN$:

\newclaim Proposition.
If $\ell_\infty$, $\ell_\infty'$ and $\ell$ are three $n$-dimensional subspaces
of $\R^{2n}$ transversal to a fourth $n$-dimensional subspace $\ell_0$, and
$\bT$ is a linear transformation from $\R^{2n}$ to itself, then
$\bN(\bT\ell_0;\bT\ell_\infty,\bT\ell) =
              \bT \bN(\ell_0;\ell_\infty,\ell)\bT^{-1}$ and
 $\bN(\ell_0;\ell_\infty',\ell) = \bN(\ell_0;\ell_\infty,\ell) +
                                              \bN(\ell_0;\ell_\infty',\ell_\infty)$.

If $\ell(t)$ is a fanning curve and $\tau$ is a real number the curve of nilpotent
operators $\bN(\ell(\tau);\ell_\infty,\ell(t))$ has a simple pole at $t = \tau$
and its residue $\bN_{-1}(\tau)$ sends $\ell_\infty$ isomorphically onto
$\ell(\tau)$. In fact, for values of $t$ close to $\tau$, $\bN$ can be
{\sl formally written\/} as the Laurent series
$$
\bN(\ell(\tau);\ell_\infty,\ell(t)) = \sum_{k=-1}^{\infty} (t-\tau)^k \bN_k(\tau),
$$
where the the coefficients $\bN_k(\tau)$ are nilpotent linear transformations
from $\R^{2n}$ to itself whose kernels contain $\ell_0$ and whose ranges are
contained in $\ell_0$. When $\ell(t)$ is spanned by a line $\cA_1 + t \cA_2$ in
the space of frames,
$$
\bN(\ell(\tau);\ell_\infty,\ell(t)) = \frac{1}{t-\tau}\bN_{-1}(\tau) + \bN_0(\tau).
$$

The notation $\bN_k(\tau)$ is somewhat inaccurate in that $\bN_{k}(\tau)$ depends on
the $(k + 2)$-jet of the curve $\ell(t)$ at $t = \tau$ and, at least a priori, on
the choice of $\ell_\infty$. However, we have the following simple result:

\newclaim Proposition.
With the previous notation, if $k \neq 0$, the transformation $\bN_k(\tau)$ does
not depend on the choice of subspace $\ell_\infty$ transversal to $\ell(\tau)$.
On the other hand, whereas $\bN_0(\tau)$ depends on this choice, the subspace
$(\bI + \bN_0(\tau))\ell_\infty$ does not.

\Proof
By Proposition~8.1, the Laurent series of $\bN((\ell(\tau);\ell_\infty',\ell(t))$
centered at  $t = \tau$ differs from that of
$\bN((\ell(\tau);\ell_\infty',\ell(t))$ only in its zeroth term
$\bN_0'(\tau) = \bN_0(\tau) + \bN(\ell_0;\ell_\infty',\ell_\infty)$.
Writing
$$
\bI + \bN_0'(\tau) = \bI + \bN_0(\tau) + \bN(\ell_0;\ell_\infty',\ell_\infty)
\quad \hbox{\rm as} \quad
(\bI + \bN_0(\tau))(\bI + \bN(\ell_0;\ell_\infty',\ell_\infty)),
$$
we see that $(\bI + \bN_0'(\tau))\ell_\infty' =  (\bI + \bN_0(\tau))\ell_\infty$.
\qed

Each one of the transformations $\bN_k(\tau)$, $k \neq 0$, is an invariant of the
fanning curve $\ell(t)$ in the sense that they only depend on its $(k+2)$-jet
at $t = \tau$, and if we change $\ell(t)$ for $\bT\ell(t)$, we change
$\bN_k$ for $\bT \bN_k(\tau) \bT^{-1}$. Using the axiomatic characterization of
the fundamental endomorphism and the horizontal curve given in Theorems~7.2 and~7.4,
it is easy to uncover the relationship between these invariants and those
considered in this paper.

\newclaim Proposition.
Let $\ell(t)$ be a fanning curve in $\grassmann$. If $\ell_\infty$ is any
$n$-dimensional subspace transversal to $\ell(\tau)$, the nilpotent transformation
$\bN_{-1}(\tau)$ coincides with the fundamental endomorphism
$\bF(\tau)$, and the subspace $(\bI + \bN_0(\tau))\ell_\infty$ is the horizontal
subspace to $\ell(t)$ at $t = \tau$.

\subsection{Cartan's method of moving frames}
While we have not found any reference that employs Cartan's method of moving
frames to study fanning curves in $\grassmann$ $(n > 1)$, it is straightforward
to generalize the approach of Flanders in \cite{Flanders:1970} to higher
dimensions.

If $\ell(t)$ is a fanning curve in $\grassmann$, we look for what Flanders calls
a {\sl natural moving frame.} This is a pair of curves of frames $\cA(t)$ and
$\cB(t)$ such that (1) $\cA(t)$ spans $\ell(t)$, (2) $\dot{\cA}(t) = \cB(t)$,
and (3) $\dot{\cB}(t) = \cA(t)R(t)$ for some curve of $n \times n$ matrices
$R(t)$. This is clearly equivalent to finding what we have called a normal
frame for $\ell(t)$. Therefore, we already know that natural moving frames always
exist and can be found at the cost of solving a linear differential equation.

The curve of $2n \times 2n$ matrices $\bA(t) = \join{\cA(t)}{\dot{\cA}(t)}$ is
what is usually called the moving frame for $\ell(t)$. By Proposition~4.2, if
$\bB(t) = \join{\cB(t)}{\dot{\cB}(t)}$ is another moving frame for $\ell(t)$,
there exists a fixed $n \times n$ invertible matrix $X$ such that
$$
\bB(t) = \bA(t) \pmatrix{X & \mzero \cr \mzero & X} .
$$

Pulling back the Maurer-Cartan form on $\lgroup{2n}$ via
the map $\bA(t)$, we obtain
$$
\bA(t)^{-1}\dot{\bA}(t) =
\pmatrix{\mzero & R(t) \cr I & \mzero}.
$$
Hence the ``invariant'' obtained Cartan's method is simply
$R(t) = -(1/2)\{\cA(t),t\}$. However, it must be taken into account that if we
had chosen the moving frame $\bB(t)$ we would have obtained $X^{-1}R(t)X$ instead.
Note that the last $n$ columns of any moving frame for $\ell(t)$
span the horizontal curve. Thus, the definition of the horizontal curve
is implicit in the moving-frame approach to the geometry of fanning curves.

\medskip
\centerline{\headfont Acknowledgements}

The authors gladly acknowledge helpful conversations with F.~Burstall,
P.~Foulon, A.~M.~Naveira, A.~Rigas, and S.~Tabachnikov. \'Alvarez-Paiva
also thanks the UNICAMP and the Universidad de Valencia for their great
hospitality during the time this paper was in preparation.

\medskip
\centerline{\headfont References}
\bibliography{../../paperbib}
\bibliographystyle{plain}

\bigskip

{\small J.~C. \'Alvarez Paiva, Department of Mathematics, Polytechnic
University, Six MetroTech Center, Brooklyn, NY, 11201.
E-mail address:} {\smalltt jalvarez@duke.poly.edu}

{\small C.~E. Dur\'an, IMECC-UNICAMP, Pra\c{c}a Sergio Buarque de Holanda,
651 Cidade Universit\'aria, Bara\~o Geraldo, Caixa Postal 6065, 13083--859
Campinas, SP, Brasil. E-mail address:} {\smalltt cduran@ime.unicamp.br}
\bye